\documentclass[12pt]{amsart}

\usepackage{amsmath,amssymb,amsfonts,amsthm,amscd}
\usepackage{manfnt}
\usepackage{mathtools}
\usepackage[mathcal]{euscript}
\usepackage{fullpage}
\usepackage{url}
\usepackage[all]{xy}
\SelectTips{cm}{10}
\usepackage{txfonts}
\usepackage[shortlabels]{enumitem}
\SetEnumerateShortLabel{C}{(C\arabic*)}
\usepackage[usenames,dvipsnames,svgnames]{xcolor}
\usepackage[colorlinks=true,linkcolor=blue,urlcolor=blue,citecolor=blue,pagebackref]{hyperref}

\renewcommand{\eqref}[1]{\textup{(\ref{#1})}}

\newtheorem{thm}{Theorem}[section]
\newtheorem{lemma}[thm]{Lemma}
\newtheorem{prop}[thm]{Proposition}
\newtheorem{cor}[thm]{Corollary}

\newtheorem{hypothesis}[thm]{Hypothesis}
\newtheorem{prop-conj}[thm]{Proposition-Conjecture}

\theoremstyle{definition}
\newtheorem{defn}[thm]{Definition}

\theoremstyle{definition}
\newtheorem{rmk}[thm]{Remark}

\theoremstyle{definition}

\theoremstyle{definition}
\newtheorem{assumption}[thm]{Assumption}

\theoremstyle{definition}

\theoremstyle{definition}

\theoremstyle{definition}

\theoremstyle{definition}

\DeclareFontFamily{U}{wncy}{}
    \DeclareFontShape{U}{wncy}{m}{n}{<->wncyr10}{}
    \DeclareSymbolFont{mcy}{U}{wncy}{m}{n}
    \DeclareMathSymbol{\Sha}{\mathord}{mcy}{"58}
\newcommand{\Q}{\mathbb{Q}}

\newcommand{\Z}{\mathbb{Z}}
\newcommand{\C}{\mathbb{C}}
\newcommand{\R}{\mathbb{R}}

\newcommand{\Qpb}{\overline{\mathbb{Q}}_p}

\DeclareMathOperator{\Hom}{Hom}

\DeclareMathOperator{\End}{End}

\DeclareMathOperator{\Sym}{Sym}

\DeclareMathOperator{\ad}{ad}
\DeclareMathOperator{\rec}{rec}
\DeclareMathOperator{\tr}{tr}
\DeclareMathOperator{\im}{im}
\DeclareMathOperator{\coker}{coker}

\DeclareMathOperator{\Spec}{Spec}

\DeclareMathOperator{\Ann}{Ann}

\DeclareMathOperator{\univ}{univ}
\DeclareMathOperator{\Fitt}{Fitt}

\newcommand{\gal}[1]{\Gamma_{#1}} % an absolute Galois group
\newcommand{\Gal}{\mathrm{Gal}} % a relative Galois group

\newcommand{\onto}{\twoheadrightarrow}

\newcommand{\mc}{\mathcal}
\newcommand{\mf}{\mathfrak}
\newcommand{\mr}{\mathrm}

\newcommand{\br}{\bar{\rho}} %to sync with lifting symplectic

\newcommand{\Frob}{\mathrm{Frob}}

\newcommand{\vpi}{\varpi}
\newcommand{\vphi}{\varphi}

\newcommand{\ov}{\overline}

\title{Effective Bounds for Adjoint Bloch--Kato Selmer Groups of Hilbert Modular Forms}

\author[C.~B.~Khare]{Chandrashekhar B. Khare}
\address{Department of Mathematics, University of California, Los Angeles,
CA 90095, USA}
\email{shekhar@math.ucla.edu}

\author[S.~Patrikis]{Stefan Patrikis}
\address{Department of Mathematics, The Ohio State University, 100 Math Tower,
231 West 18th Avenue, Columbus, OH 43210, USA}
\email{patrikis.1@osu.edu}


\begin{document}

\begin{abstract}
Let $F$ be a totally real field and let $\Pi$ be a regular algebraic cuspidal
automorphic representation of $\mathrm{GL}_2(\mathbb A_F)$.  Under the usual
Taylor--Wiles residual irreducibility hypothesis, we prove the vanishing of the
adjoint Bloch--Kato Selmer group with characteristic-zero coefficients and give
an explicit upper bound for a corresponding Selmer group with divisible
coefficients.  The argument revisits Wiles's use of level raising at
Taylor--Wiles primes.  In place of freeness or almost freeness of automorphic
cohomology over a Hecke algebra, we compare congruence modules at the original
and auxiliary levels and control the part of the Hecke cotangent space detected
by the automorphic module.  This produces Galois cohomology classes with
prescribed ramification at the auxiliary primes.  The resulting bounds depend
explicitly on congruence data at the original automorphic level and on torsion
in cotangent spaces of local deformation rings.
\end{abstract}

\maketitle

\section{Introduction}
Let $\Pi$ be a regular C-algebraic cuspidal automorphic representation of $\mr{GL}_2(\mathbb{A}_F)$ over a totally real field $F$, and let $\rho= \rho_{\Pi, \iota} \colon \Gal(\overline{F}/F) \to \mr{GL}_2(\Qpb)$ be the $p$-adic Galois representation associated to $\Pi$ (for a choice of isomorphism $\iota \colon \Qpb \to \C$). For some finite extension $E/\Q_p$, with ring of integers $\mc{O}$, we can assume (replacing $\rho$ by a conjugate) $\rho$ factors through $\mr{GL}_2(\mc{O})$. We can then consider different versions of the Bloch-Kato Selmer groups for the trace-zero adjoint representation. Let $W= \ad^0(\rho)$, $W_E= \ad^0(\rho) \otimes E$, and $W_{E/\mc{O}}= \ad^0(\rho) \otimes_{\mc{O}} E/\mc{O}$. Then if $S$ is a finite set of places such that $\Pi$ is unramified outside $S$, and we let $\gal{F, S}= \Gal(F(S)/F)$ denote the Galois group of the maximal unramified outside $S$ extension of $F$, we first have the Selmer groups $H^1_g(\gal{F, S}, W_E)$ and $H^1_f(F, W_E)$: the former classifies lifts of $\rho \otimes E$ to $\mr{GL}_2(E[\epsilon])$, with fixed determinant $\det(\rho)$, that are unramified outside $S$ and de Rham at all places $v \vert p$; the latter consists of cohomology classes that are unramified away from $p$ and at $v \vert p$ lie in the local Bloch-Kato group $H^1_f (\gal{F_v}, W_E)$, which in fact in this case coincides with the previous local condition $H^1_g(\gal{F_v}, W_E)$ that the lift be de Rham (see \cite[Remark 1.2.9]{allen:selmer}). Then the Bloch-Kato conjecture (building on conjectures of Deligne, Beilinson, and Fontaine-Mazur) in this setting first asserts that the groups $H^1_f(F, W_E)$ and $H^1_g(\gal{F, S}, W_E)$ vanish, a prediction whose automorphic interpretation is the familiar fact that $L(\ad^0(\Pi), 1)$ is non-zero.

Bloch and Kato's innovation was moreover, at least when one has a motivic interpretation of $\rho$, to predict the special value of the L-function in terms of a transcendental period, local Tamagawa factors (arising from places in $S$ where $\rho$ is ramified), and the torsion variant of the above Selmer group, $H^1_f(F, W_{E/\mc{O}})$, where the local conditions are by definition the images of the local conditions for $W_E$; this can be translated to predictions about the size of the more relaxed Selmer groups $H^1_g(\gal{F, S}, W_{E/\mc{O}})$ via global duality (as in \cite[Lemma 2.1]{diamond-flach-guo}). Galois deformation theory relates both the rational and the torsion Selmer groups to tangent spaces of suitable Galois deformation rings $R$ for $\br$ with local conditions at the primes in $S$. Set up properly, if $\phi \colon R \to \mc{O}$ is the augmentation corresponding to the deformation $\rho$, with $\mf{p}= \ker(\phi)$, the tangent space $\Hom(\mf{p}/\mf{p}^2, E)$ will be isomorphic to $H^1_g(\gal{F, S}, W_E)$, and its torsion variant $\Hom(\mf{p}/\mf{p}^2, E/\mc{O})$ will be some quantifiable variant on $H^1_g(\gal{F, S}, W_{E/\mc{O}})$.

Thus we are faced with a hierarchy of possible results toward the Bloch--Kato conjecture in this or related settings:
\begin{enumerate}
\item Prove the vanishing of $H^1_g(\gal{F, S}, W_E)$, or equivalently the finiteness of $H^1_g(\gal{F, S}, W_{E/\mc{O}})$.
\item Prove $H^1_g(\gal{F, S}, W_{E/\mc{O}})$ or $\Hom(\mf{p}/\mf{p}^2, E/\mc{O})$ is finite, with an explicit bound  or even exact formula for the order.
\item Prove the full Bloch-Kato conjecture for $H^1_f(F, W_{E/\mc{O}})$.
\end{enumerate}
In his work on Fermat's Last Theorem, Wiles proved many cases of (1), along with (in combination with the main theorem of \cite{taylor-wiles:fermat}) an exact formula for corresponding cases of (2) (see \cite[Theorem 3.3]{wiles:fermat}). Wiles's work depends on a number of simplifying hypotheses: the field $F$ is $\Q$, the residual representation $\br$ has ``large" image as required by the traditional Taylor--Wiles method, and the local representations $\rho|_{\gal{F_v}}$ satisfy simplifying hypotheses, particularly at $v \vert p$ (e.g., $\Pi$ should be weight 2 and either unramified at $p$ or ordinary at $p$). The structure of Wiles's formula depends on proving an $R=\mathbb{T}$ theorem identifying $R$ with an appropriate reduced Hecke algebra $\mathbb{T}$, and on the fact that these rings are complete intersections. This allows, via Wiles's numerical criterion, the size of $\mf{p}/\mf{p}^2$ (now known to be finite, since $\mathbb{T}$ is reduced) to be reinterpreted as the size of $\mc{O}/\eta_{\mathbb{T}}$, where $\eta_{\mathbb{T}}$ is the congruence ideal. 

The mechanism of Wiles's proof is to introduce auxiliary (Taylor-Wiles) primes $Q$ of ramification in the deformation problem and show using automorphic congruences that the growth of the Selmer group with relaxed ramification at $Q$ is maximal. In \cite[Theorem~3.1(i)]{wiles:fermat}, this method is used to prove the upper bound in (2); the full force of an integral $R=\mathbb{T}$ theorem, relying on the Taylor-Wiles patching method, can then be used to deduce that the upper bound is exact (\cite[Theorem~3.3]{wiles:fermat}).  While Wiles does not exactly establish (3), Diamond--Flach--Guo (\cite{diamond-flach-guo}) were able to prove (3) for $F=\Q$, under the hypotheses $p>k$ and $\Pi_p$ unramified, by analyzing an explicit motivic setup in combination with the Taylor-Wiles method and an idea of Hida that interprets the special $L$-value in terms of a congruence ideal.

We describe some of the subsequent work on questions (1)-(3). Kisin (\cite{kisin:geometricdef}), still when $F= \Q$, took up Wiles's method and proved the vanishing of $H^1_g(\gal{\Q, S}, \ad(\rho))$ for $\Pi$ of arbitrary weight and power of $p$ in the level. His proof requires a difficult analysis of the integral $p$-adic Hodge theory of modular curves but is very close in spirit to Wiles's method. Kisin's method does not yield an explicit bound as in (2), but it follows more closely the argument of \cite[Chapter 3]{wiles:fermat} and in particular does not rely on an $R=\mathbb{T}$ theorem of any kind. Later works have relied heavily on Taylor-Wiles patching and have proven only (1). We highlight the work of Allen (\cite{allen:selmer}), which proves vanishing of $H^1_g(\gal{F, S}, W_E)$ in the much more general setting where $W= \ad(\rho_{\Pi, \iota})$ for $\Pi$ a regular algebraic polarizable (in the sense of \cite{blggt:potaut}) cuspidal automorphic representation of $\mr{GL}_n(\mathbb{A}_K)$ where $K$ is a CM field with totally real subfield $F$;\footnote{Note that $\rho_{\Pi, \iota}$ is defined on $\gal{K}$, but there is an outer action of $\gal{F}$ involved in the definition of the Selmer group, which only classifies suitably polarized deformation.} Allen's work allows for arbitrary local behavior of $\Pi$, but it requires $p > 2$ and that $\br(\gal{K(\zeta_p)})$ have large (adequate) image to allow the usual higher-rank Taylor-Wiles method (developed in \cite{clozel-harris-taylor}, \cite{thorne:adequate}) to be applied. In a major innovation, Newton and Thorne (\cite{newton-thorne:selmer}) and Thorne (\cite{thorne:selmer}) have proven (1) for all regular algebraic polarizable cuspidal automorphic representations of $\mr{GL}_n(\mathbb{A}_K)$, subject only to the very weak requirement that $\rho_{\Pi, \iota}|_{\gal{K(\zeta_{p^\infty})}}$ be irreducible (their work has recently been extended to the non-polarized case in \cite{acampo-hevesi-thorne-whitmore}). Results along the lines of (2) (let alone (3)) seem to be inaccessible via these methods. Forthcoming work of Skinner-Sweeting finds sharp bounds in parallel weight 2. In the direction of (3), work in progress of Diamond-Iyengar-Khare-Manning seeks exact Bloch-Kato formulas beyond the setting of \cite{diamond-flach-guo}, but still restricted to $F= \Q$, using integral $R=\mathbb{T}$ theorems and the $p$-adic local Langlands correspondence for $\mr{GL}_2(\Q_p)$.

Our work, which in the present paper is restricted to $\mr{GL}_2$ over our totally real field $F$, revisits the arguments of \cite[Chapter 3]{wiles:fermat} and gives a framework where we can prove explicit upper bounds for $H^1_g(\gal{F, S}, W_{E/\mc{O}})$ as in (2) for regular algebraic cuspidal automorphic representations $\Pi$ of ${\mathrm GL}_2({\mathbb A}_F)$.  Our method, like the argument of Wiles in \cite[Chapter 3]{wiles:fermat}, is not able to give formulas for the exact order of this Selmer group. We also rely on the Taylor--Wiles hypothesis that $\br|_{\gal{F(\zeta_p)}}$ is absolutely irreducible, but are otherwise able to prove meaningful upper bounds with no restriction on the local behavior of $\Pi$. See \S \ref{generalsection} for possible generalizations: we are optimistic that via these methods we can establish explicit upper bounds as in (2) but in the generality of Allen's work, and we have interesting lines of inquiry that may eventually make possible a version of (2) in the Newton--Thorne setting.

More precisely, we work on a quaternion algebra $D$ that is ramified precisely at infinity, and with automorphic representations of the associated algebraic group $G^{D}$ over $F$ (see \S \ref{generalsetup}); in general we reduce to this case by base-change, which is harmless for establishing (1) but leads to some inefficiency in (2). We refer to Theorem \ref{finalthm} for the result in the Hilbert modular case, and for the rest of the introduction we assume that $\Pi$ arises from an automorphic representation $\pi$ of $G^D(\mathbb{A}_F)$ via the Jacquet-Langlands correspondence. Letting $U \subset G^D(\mathbb{A}_F^\infty) \cong \mr{GL}_2(\mathbb{A}_F^\infty)$ be the conductor subgroup, we have a space of $p$-adically integral automorphic forms $S(U, \mc{O})$ on $G^D$ to which $(\pi^\infty)^U$ contributes, and we have the reduced Hecke algebra $\mathbb{T}(U)$ generated by good Hecke operators $T_v$, $v \not \in S$ acting on $S(U, \mc{O})$. Since $(\pi^\infty)^U$ is 1-dimensional, we obtain an eigencharacter $\lambda \colon \mathbb{T}(U) \to \mc{O}$ giving the Hecke eigenvalues. Let $I= \ker(\lambda)$, and let $\mf{m} = \ker(\lambda \pmod{\vpi})$. Our results will be stated in terms of $M_{\emptyset}= S(U, \mc{O})_{\mf{m}}$ and $I \subset \mathbb{T}_{\emptyset}=\mathbb{T}(U)_{\mf{m}}$. In \S 3 we will define Selmer conditions $\mc{D}_S$ at all torsion levels and the resulting torsion Selmer group $H^1_{\mc{D}_S}(F, W_{E/\mc{O}})$ (see the discussion following Definition \ref{defring}). For a finite $\mc{O}$-module $A$, write $\lg(A)=\operatorname{length}_{\mc{O}}(A)$; when a cohomology group is finite, write $h^i_{\mc{L}}(F,V)=\lg(H^i_{\mc{L}}(F,V))$. Then we show:
\begin{thm}\label{introthm}
Assume that $[F(\zeta_p):F] > 2$, that $\br|_{\gal{F(\zeta_p)}}$ is absolutely irreducible, and if $p=5$ and $\br$ has projective image $\mr{PGL}_2(\mathbb{F}_5)$, then $F(\zeta_5)$ is not contained in $\ov{F}^{\ker(\mr{proj}(\br))}$.\footnote{We can allow $p=3$ if $U$ is assumed sufficiently small; the method would let us take $p=3$ with no restriction on $U$, but the final bound would be different, and the technical complications obfuscating for our present purposes.} Then $H^1_g(\gal{F, S}, W_E)=0$, and $H^1_{\mc{D}_S}(\gal{F, S}, W_{E/\mc{O}})$ is finite with length (as $\mc{O}$-module) at most $c_S+d_r$ where:
\begin{itemize}
\item $c_S= \sum_{v \in S} c_v$ is a sum of local error terms arising from singularities of local deformation rings: $c_v= \lg(\Omega_{\ov{R}_v^{\square}/\mc{O}} \otimes_{\ov{R}^{\square}_v, \rho} \mc{O})_{\mr{tor}}$ where for $v \nmid p$, $\ov{R}^{\square}_v$ is the local lifting ring with fixed determinant $\det(\rho)$ for $\br|_{\gal{F_v}}$; and for $v \vert p$ it is a certain potentially semistable deformation ring specified in Lemma \ref{gw}. 
\item $r=\dim_k H^1_{\mc{D}_S^\perp}(\gal{F, S}, \ad^0(\br)^*) +\sum_{v \in S} t_v$; here $t_v$ (defined in Equation \eqref{localtorsion}) is the number of generators of $(\Omega_{\ov{R}_v^{\square}/\mc{O}} \otimes_{\ov{R}^{\square}_v, \rho} \mc{O})_{\mr{tor}}$.
\item $d_r= r(j_{\emptyset}+a_{\emptyset})+b_{\emptyset}- \lg\bigl(\Psi_{\lambda}(M_{\emptyset})\bigr)$, where
\begin{itemize}
\item $a_{\emptyset}= \lg\bigl((M_{\emptyset}/IM_{\emptyset})_{\mr{tor}}\bigr)$.
\item $b_{\emptyset}= \lg(IM_{\emptyset}/I^2M_{\emptyset})$, and
\[
\vpi^{j_{\emptyset}}\mc{O}
=\Ann(IM_{\emptyset}/I^2M_{\emptyset})
 \Ann((M_{\emptyset}/IM_{\emptyset})_{\mr{tor}})^2.
\]
\item $\Psi_{\lambda}(M_{\emptyset})$ is the congruence module defined in \eqref{congruence-module}.
\end{itemize}
\end{itemize}
\end{thm}
\begin{rmk}
Let us demystify this formula in the simple case where $M_{\emptyset}$ is free of rank 1 over $\mathbb{T}_{\emptyset}$ (analogous to Wiles's situation where his Hecke algebra is Gorenstein, and the localized cohomology of the modular curve is free rank 2 over the Hecke algebra; in our definite quaternion algebra setting, the generic ranks will be 1 rather than 2). %Here $R_{\emptyset}$ would be a deformation ring that is a quotient (with perhaps somewhat stricter local conditions) of the one implicit in the preceding discussion (and explicit in Definition \ref{defring} below). 
Then $a_{\emptyset}=0$, $\Psi_{\lambda}(M_{\emptyset})\simeq \mc{O}/\eta_{\lambda}(\mathbb{T}_{\emptyset})$, and $b_{\emptyset}= \lg(I/I^2)$. %\lg(\mf{p}_{\emptyset}/\mf{p}_{\emptyset}^2)$. 
Theorem \ref{introthm}, %or rather its proof with $R_{\emptyset}$ in place of the ring we use in Definition \ref{defring}, 
would then bound the length of $H^1_{\mc{D}_S}(\gal{F, S}, W_{E/\mc{O}})$ by
\[
\dim_k H^1_{\mc{D}_S^\perp}(\gal{F, S}, \ad^0(\br)^*) j_{\emptyset}+ \lg(I/I^2)- \lg\bigl(\mc{O}/\eta_{\lambda}(\mathbb{T}_{\emptyset})\bigr)+c_S. 
\]
If we further simplify to the case where the $c_v$ (and thus $t_v$) are all zero, then the Greenberg-Wiles formula transforms this bound into
\[
\dim_k H^1_{\mc{D}_S}(\gal{F, S}, \ad^0(\br)) \cdot j_{\emptyset}+ \lg(I/I^2)- \lg\bigl(\mc{O}/\eta_{\lambda}(\mathbb{T}_{\emptyset})\bigr).
\]
This mod $p$ Selmer dimension gives the maximum number of generators needed for the deformation ring and thus bounds from above the number of $\mc{O}$-module generators of $I/I^2$, so this $\dim_k H^1_{\mc{D}_S}(\gal{F, S}, \ad^0(\br)) \cdot j_{\emptyset}$ term is a weakened version (expected number of generators times exponent being an upper bound on length) of $\lg(I/I^2)$. Wiles's sharp formula in \cite[Theorem 3.1]{wiles:fermat} is $2 \lg(I/I^2)- \lg\bigl(\mc{O}/\eta_{\lambda}(\mathbb{T}_{\emptyset})\bigr)$. Thus our general formula is not quite as strong as Wiles's, but it is conceptually closely related. It would be interesting to expand on this connection and try to refine our arguments to yield a sharper final result in certain simplified settings.
\end{rmk}
Our method, like Wiles's, produces, for each $n \geq 1$, Taylor-Wiles sets $Q_n$ of level $n$ and enough ramified Galois cohomology classes in $H^1_{\mc{D}_{S \cup Q_n}}(\gal{F, S \cup Q_n}, \ad^0(\rho) \otimes \mc{O}/p^n)$ (the local conditions as before in $S$ but now allowing all local deformations at $Q_n$) to (provided $Q_n$ is chosen to annihilate the dual Selmer group mod $p^n$) bound $H^1_{\mc{D}_S}(\gal{F, S}, \ad^0(\rho) \otimes \mc{O}/p^n)$ using the Greenberg-Wiles formula (see Lemma \ref{gw} and Proposition \ref{mainbound}). At each Taylor-Wiles level $U_{Q_n}$, we have a reduced Hecke algebra with diamond operators $\mc{O}[\Delta_{Q_n}] \to \mathbb{T}_{Q_n}$ acting on a space $M_{Q_n}$ of automorphic forms, and we have a choice of stabilization $\lambda_{Q_n} \colon \mathbb{T}_{Q_n} \to \mc{O}$ of $\lambda$.  The linchpin of the method is Proposition \ref{twsetup} and Proposition \ref{growth}, which show in a rather direct manner that $\mf{a}_{Q_n}/\mf{a}_{Q_n}^2 \to I_{Q_n}/I_{Q_n}^2$ ($I_{Q_n}= \ker(\lambda_{Q_n})$) has large image by considering their respective images in $\Hom(M_{Q_n}/I_{Q_n}M_{Q_n}, I_{Q_n}M_{Q_n}/I^2_{Q_n} M_{Q_n}) \subset \End(M_{Q_n}/I^2_{Q_n}M_{Q_n})$. The method depends in an essential way on $M_{Q_n}/I_{Q_n}M_{Q_n}$ having generic rank 1 and on the length of its torsion submodule being bounded independently of $Q_n$; in the Taylor-Wiles setting used here $M_{Q_n}/I_{Q_n}M_{Q_n}$ is in fact independent of $Q_n$.

%We close this introduction by noting that the idea of bounding arithmetic invariants by constructing ramified cohomology classes at auxiliary primes has an august history going back at least to the work of Iwasawa, who showed that Leopoldt's conjecture is equivalent to .

We normalize class field theory to send uniformizers to geometric Frobenius elements. We normalize Hodge-Tate weights so that the cyclotomic character $\kappa$ has Hodge-Tate weight 1.
\section{Setup and preliminary reductions}\label{generalsetup}
Let $F_0$ be a totally real field, and let $\Pi_0$ be a regular C-algebraic cuspidal automorphic representation of $\mr{GL}_2(\mathbb{A}_{F_0})$. Fix an isomorphism $\iota \colon \ov{\Q}_p \xrightarrow{\sim} \C$. Via $\iota$, each embedding $\tau \colon F_0 \to \ov{\Q}_p$ induces $\iota \tau$ and thus a unique $v \vert \infty$; since $F_0$ is totally real, there is no ambiguity in writing $\tau_v \colon F_0 \to \R$ for $\iota \tau$. The Harish--Chandra parameter of $\Pi_{0,v}$ has the form $\{(-m_{\tau_v}+\frac{1}{2},-n_{\tau_v}-m_{\tau_v}-\frac{1}{2})\}_{v\mid\infty}$ for some integers $n_{\tau_v}\geq0$ and $m_{\tau_v}$ (this normalization will be convenient below). We write $n_{\tau}=n_{\tau_v}$ and $m_{\tau}=m_{\tau_v}$, where $\tau$ and $\tau_v$ are associated as above.

By the work of many people, notably Carayol, Taylor, Tsuji, Saito, and Skinner (see \cite{skinner:hmf} for a summary), there is a continuous representation $\rho_{\Pi_0, \iota} \colon \gal{F_0} \to \mr{GL}_2(\ov{\Q}_p)$ such that 
\begin{enumerate}
\item $\rho_{\Pi_0, \iota}$ is odd.
\item For all finite places $v \nmid p$, $\mr{WD}(\rho_{\Pi_0, \iota}|_{\gal{F_{0,v}}})^{F\text{-}\mathrm{ss}} \cong \iota(\rec_v(\Pi_{0,v}|\det|^{-1/2}))$.
\item For $v \vert p$, $\rho_{\Pi_0, \iota}|_{\gal{F_{0, v}}}$ is de Rham, and the $\tau$-labeled Hodge-Tate weights are $\{-m_{\tau}, -n_{\tau}-m_{\tau}-1 \}$. Moreover, this representation is semistable (respectively, crystalline) when $\Pi_{0, v}$ is an unramified twist of the Steinberg representation (respectively, an unramified representation).
\end{enumerate}
We fix an integral model, still denoted $\rho_{\Pi_0, \iota} \colon \gal{F_0} \to \mr{GL}_2(\mc{O})$, defined over the ring of integers $\mc{O}$ in a finite extension $E/\Q_p$, which we assume large enough to contain the images of all embeddings $F_0\hookrightarrow\overline{\Q}_p$. After choosing the auxiliary field $F$ below, we enlarge $E$ further if necessary. At one point later in the argument (part (5) of Proposition \ref{twsetup}), if $p$ is less than the difference of the Hodge-Tate weights at some embedding, we will need to replace $\mc{O}$ by its quadratic extension $\mc{O}[\sqrt{\vpi}]$; we do this from now, so that we don't have to change $\mc{O}$ later on. We continue to denote the resulting coefficient field and its ring of integers by $E$ and $\mc O$, and let $e$ be the ramification index of this coefficient field over $\Q_p$.

Let $S_0$ be the finite set of places of $F_0$ consisting of the infinite places, the
places above $p$, and the places at which $\Pi_0$ is ramified.  For any
extension $F/F_0$ considered below, let $S$ denote the set of places of
$F$ above $S_0$.

\iffalse
There is a finite solvable extension $F/F_0$ such that:
\begin{enumerate}
\item $[F: \Q]$ is even.
\item For all finite places $v$ of $F$, the base change $\mr{BC}_{F/F_0}(\Pi_0)_v$ is either unramified or is an unramified twist of the Steinberg representation. 
\item
\end{enumerate}
\fi
Let $F/F_0$ be any totally real finite solvable extension such that $[F:\Q]$ is even (later on we will make more requirements on $F$). Let $\Pi= \mr{BC}_{F/F_0}(\Pi_0)$. It similarly has an associated Galois representation $\rho_{\Pi, \iota}$, and we will mainly analyze this restriction. Let $D$ be the quaternion algebra over $F$ whose set of ramified places is exactly $S_{\infty}= \{ v \vert \infty\}$. Let $G=G^{D}$ be the reductive group over $F$ associated to $D$, so that $G^D(R)=(D\otimes_F R)^\times$ for every $F$-algebra $R$. Fix a maximal order $\mc{O}_D$ of $D$ along with isomorphisms $\mc{O}_{D, v} \xrightarrow{\sim} M_2(\mc{O}_{F_v})$ for all finite places $v$; from now on we make this identification implicitly and use it to identify $G(\mathbb{A}_F^\infty) \cong \mr{GL}_2(\mathbb{A}_F^\infty)$. Then $\Pi$ lies in the image of the Jacquet--Langlands correspondence from infinite-dimensional cuspidal automorphic representations of $G(\mathbb{A}_F)$ to cuspidal automorphic representations of $\mr{GL}_2(\mathbb{A}_F)$, and we let $\pi$ be the cuspidal automorphic representation of $G(\mathbb{A}_F)$ such that $\mr{JL}(\pi)= \Pi$. Let $\chi$ be the common central character of $\pi$ and $\Pi$. Let $U= \prod_v U_v \subset G(\mathbb{A}_F^\infty)$ be the (compact open) conductor subgroup of $\pi$, so $(\pi^\infty)^U$ is one-dimensional.

%At the infinite places, there are integers $\{k_{\tau}, s_{\tau}\}_{\tau \colon F \to \R}$ such that 
By the explicit form of the Jacquet-Langlands correspondence, for each $v \vert \infty$, $\pi_v \cong \Sym^{n_{\tau_v}} (\C^2) \otimes \det^{-n_{\tau_v}-m_{\tau_v}}$. As notational shorthand, we will write $W_{\iota}= \otimes_{v \vert \infty} \Sym^{n_{\tau_v}}(\C^2) \otimes \det^{m_{\tau_v}}$, isomorphic to $\pi_{\infty}^\vee$. Necessarily $w= n_{\tau}+2 m_{\tau}+1$ is independent of $\tau$, and note then that $\chi|_{F_\infty^\times}=N_{F_{\infty}/\R}^{1-w}$. Associated with $\chi$ we have the $p$-adic character $\chi_{\iota} \colon \mathbb{A}^\times_F/\ov{F^\times F_{\infty}^\times} \to \ov{\Q}_p^\times$ given by
\[
\chi_{\iota}(x)= \prod_{\tau \colon F \to \ov{\Q}_p} \tau(x_{v_{\tau}})^{1-w} \cdot \iota^{-1} \left( N_{F_{\infty}/\R}(x_{\infty})^{w-1} \cdot \chi(x) \right),
\]
and by the global Artin map (normalized to send geometric Frobenii to uniformizers) we may regard $\chi_\iota$ as a Galois character $\gal{F} \to \ov{\Q}_p^\times$. We then see that $\det(\rho_{\Pi, \iota})= \kappa^{-1} \chi_{\iota}$, where $\kappa$ denotes the $p$-adic cyclotomic character.

We can now introduce the spaces of ($p$-adic) integral automorphic forms. A slight technicality intervenes, since we will want to have $\mc{O}$-linear self-dualities on these spaces, even when some $n_{\tau}$ is bigger than $p-1$. Thus following \cite[\S 3.1.9]{kisin:bt} and \cite[Lemma 3.3.1]{conrad-diamond-taylor:modularity}, we note that for any $n \geq 0$ there is an $\mc{O}$-lattice $L_n \subset \Sym^n(E^2)$ that is stable under $M_2(\mc{O})$ and on which the usual $E$-linear duality $\Sym^n E^2 \times \Sym^n E^2 \to \det^{n}$ restricts to an $\mc{O}$-linear perfect pairing. For any $\mc{O}$-algebra $A$, let $L_{n, A}= L_n \otimes_{\mc{O}} A$, and consider the $A$-module $S(A)= S_{\underline{n}, \chi_{\iota}}(A)$ of functions
\[
\varphi \colon D^\times \backslash (D \otimes \mathbb{A}_F^{\infty})^\times \to \otimes_{ \tau \colon F \to \ov{\Q}_p} L_{n_{\tau}, A} \otimes {\det}^{m_{\tau}}(A^2) %previously \Sym^{n_{\tau}}(A^2)
\]
such that
\begin{itemize}
\item $\vphi(gz)=\chi_{\iota}(z) \vphi(g)$ for all $g \in (D \otimes \mathbb{A}_F^\infty)^\times$ and $z \in (\mathbb{A}_F^\infty)^\times$; and
\item there exists a compact open subgroup $V \subset G(\mathbb{A}_F^\infty)$ such that $\vphi(gv)= v_p^{-1} \vphi(g)$ for all $v \in V$, $g \in G(\mathbb{A}_F^\infty)$.
\end{itemize}
When $A= \mc{O}$, we let $\langle \cdot, \cdot \rangle$ be the perfect pairing induced on $\otimes_{ \tau \colon F \to \ov{\Q}_p} L_{n_{\tau}} \otimes {\det}^{m_{\tau}}$ by the above pairings; we will only need this in Proposition \ref{twsetup}. We define an action of $G(\mathbb{A}_F^{\infty, p}) \times (\mr{GL}_2(F_p) \cap M_2(\mc{O}_{F, p}))$ on $S(A)$ by $(u \cdot \vphi)(g)= u_p \vphi(gu)$, where $u_p$ denotes projection to the $p$-component. When $A$ is an $E$-vector space, this extends to an action of $G(\mathbb{A}_F^\infty)$. For a particular compact open subgroup $V$, we write $S(V, A)$ for the $V$-invariants with respect to this action. Then we have the Hecke operators $[VgV] \in \End_A(S(V, A))$, given by $[VgV]\vphi= \sum x_i \cdot \vphi$ where $VgV= \bigsqcup_i x_i V$, which make sense in general for $g \in G(\mathbb{A}_F^{\infty, p}) \times (\mr{GL}_2(F_p) \cap M_2(\mc{O}_{F, p}))$, and when $A$ is an $E$-vector space make sense for any $g \in G(\mathbb{A}_F^\infty)$.

There is a $G(\mathbb{A}_F^\infty)$-equivariant isomorphism
\begin{equation}\label{cusp}
j \colon S(\mc{O}) \otimes_{\mc{O}, \iota} \C \xrightarrow{\sim} \Hom_{G(F_{\infty})}(W_{\iota}^\vee, \mc{A}_{G, \chi})\xrightarrow{\sim} \bigoplus_{\sigma} \sigma^\infty,
\end{equation}
where $\mc{A}_{G, \chi}$ denotes the space of (complex) automorphic forms on $G(\mathbb{A}_F)$ with central character $\chi$, and where the direct sum is taken over automorphic representations $\sigma$ of $G(\mathbb{A}_F)$ with central character $\chi$ and such that $\sigma_{\infty}$ is isomorphic to $\pi_{\infty}$.\footnote{Define $j(\varphi)(w^\vee)(g)= w^\vee(g_{\infty}^{-1} g_p\vphi(g^\infty))$.} Taking $U$-invariants, we obtain an analogous isomorphism, now equivariant for the Hecke operators. Let $\mathbb{T}(U)$ be the reduced Hecke algebra at level $U$, defined as the $\mc{O}$-subalgebra of $\End(S(U, \mc{O}))$ generated by all $T_v$, $v \not \in S$.

Now, by the choice of $U$, $(\pi^\infty)^U$ is 1-dimensional, so we obtain the eigencharacter $\lambda \colon \mathbb{T}(U) \to \mc{O}$ such that $\lambda(T_v)= N(v)^{1/2}(\chi_{1, v}(\vpi_v)+\chi_{2, v}(\vpi_v))$ when $\pi_v \cong \mr{ind}_B^{G(F_v)} (\chi_{1, v} \times \chi_{2, v})$ (normalized induction).
%, and $\lambda(U_{\vpi_v})= N(v) \chi_v(\vpi_v)$ when $\pi_v \cong \mr{St} \otimes (\chi_v \circ \det)$. 
Reducing mod $\vpi$, we obtain a maximal ideal $\mf{m} \subset \mathbb{T}(U)$, and we let $\mathbb{T}= \mathbb{T}(U)_{\mf{m}}$ be the localization, with $\lambda \colon \mathbb{T} \to \mc{O}$ the homomorphism induced by $\lambda$, and with $I= \ker(\lambda)$.

\iffalse
\begin{lemma}
There is a group determinant $D_{\mf{m}} \colon \gal{F, S} \to \mathbb{T}$ characterized uniquely by $D_{\mf{m}}(\Frob_v)=T_v$ for all $v \not \in S$. When $\mathfrak{m}$ is non-Eisenstein, there is a continuous representation $\rho_{\mf{m}} \colon \gal{F, S} \to \mr{GL}_2(\mathbb{T})$ such that $\tr(\rho_{\mf{m}})= D_{\mf{m}}$.
\end{lemma}
Let $\ov{D}= D_{\mf{m}} \pmod {\mf{m}}$ be the reduction, and let $R_{\ov{D}, S}$ be the complete local Noetherian $\mc{O}$-algebra representing the functor which associates to $A \in \mc{C}_{\mc{O}}$ the set of continuous group determinants $D$ of $\gal{F, S}$ over $A$ such that
\begin{itemize}
\item $D \pmod {\mf{A}}= \ov{D}$.
\item For all $g \in \gal{F, S}$, $D(t-g)$ has constant term $\kappa^{-1}\chi_{\iota}(g)$. 
\end{itemize}
Let $a= \min_{\tau}()$, and let $b= \max_{\tau}()$. Then we write $R^{[a, b]}_{\ov{D}, S}$ for the quotient representing those $D$ that further satisfy:
\begin{itemize}
\item There is a Cayley-Hamilton representation $(\mc{O}[\gal{F, S}], D) \to (B, D')$ over $A$ such that for each $n \geq 1$, $B/\mf{m}_A^n B$ lies in $\mc{E}_{F, S}^{[a, b]}$, the category of finite $\Z_p[\gal{F, S}]$-modules whose restrictions to $\gal{F_v}$ for all $v \vert p$ are isomorphic to quotients of lattices in semistable representations with Hodge-Tate weights in the interval $[a, b]$.
\end{itemize}
\fi

We will also need these spaces with level structure added at auxiliary (Taylor-Wiles) primes. First let $Q$ be any set of primes of $F$ disjoint from $S$. We first introduce the Taylor-Wiles level subgroups 
\[
U_{Q, v}= \left\{\begin{pmatrix} a & b \\ c & d \end{pmatrix} \in U_0(v): \frac{a}{d} \mapsto 1 \in \Delta_v \right\}, 
\]
where $\Delta_v$ is the $p$-Sylow subgroup (equivalently, quotient) of the cyclic group $\kappa(v)^\times$, and set $U_Q= \prod_{v \not \in Q} U_v \times \prod_{v \in Q} U_{Q, v}$. $S(U_Q, \mc{O})$ has commuting actions of $\Delta_Q= \prod_{v \in Q} \Delta_v$ and $\mathbb{T}_Q(U_Q)$, the commutative $\mc{O}$-subalgebra of $\End(S(U_Q, \mc{O}))$ generated by all $T_v$, $v \not \in S \cup Q$ and by $U_{\vpi_v}$, $v \in Q$, for some fixed choice of uniformizers $\vpi_v$ of $\mc{O}_{F_v}$.

Provided $U_0(Q)$ is sufficiently small, $S(U_Q, \mc{O})$ is a finite free $\mc{O}[\Delta_Q]$-module by \cite[Prop. 5.4]{gee:aws}, and $S(U_Q, \mc{O})/\Delta_Q \xrightarrow{\sim} S(U_0(Q), \mc{O})$ (with the map induced by summing over the $\Delta_Q$-action). In \S \ref{twsection}, we will work under an Assumption (Assumption \ref{twassumption}) that guarantees this always holds.

Choose a root $\alpha_v$ of $X^2- \lambda(T_v) X + N(v) \chi(\vpi_v)$. For any $v \in Q$, $\pi_v^{U_{Q, v}}$ is two-dimensional, and the $\alpha_v$-eigenspace of $U_{\vpi_v}$ is one-dimensional (whether or not the roots of the Hecke polynomial are distinct). %in a suitable basis by the matrix in the proof of Lemma \ref{stronggenfree}, so it has a \textit{unique} (whether or not the Hecke polynomial has distinct roots) eigenvector in $\pi_v^{U_{Q, v}}$ with eigenvalue $\alpha_v$. 
We obtain a character $\lambda_Q \colon \mathbb{T}_Q(U_Q) \to \mc{O}$ 
%(here $\mathbb{T}(U_Q)$ is defined just as $\mathbb{T}(U)$ was, with the level subgroup $U_Q$ in place of $U$ in the previous definition) 
corresponding to the choices $(\alpha_v)_{v \in Q}$ (and the initial choice of $\lambda$). Let $I_Q= \ker(\lambda_Q)$, and let $\mf{m}_Q = \ker(\lambda_Q \pmod {\vpi})$.

Define $M_Q= S(U_Q, \mc{O})_{\mf{m}_Q}$; it is a finite $\mathbb{T}_Q:= \mathbb{T}_Q(U_Q)_{\mf{m}_Q}$-module, finite free as an $\mc{O}$-module.
\begin{lemma}\label{weakgenfree}
$(M_Q/I_QM_Q) \otimes_{\mc{O}} E$ is one-dimensional over $E$.
\end{lemma}
\begin{proof}
The decomposition of Equation \eqref{cusp} shows that  $S(U_Q, E) \otimes_{\iota} \C$ is $\mathbb{T}_Q(U_Q)$-equivariantly isomorphic to $\bigoplus_{\sigma} \sigma^{U_Q}$, where the sum is indexed over isomorphism classes of irreducible admissible representations $\sigma$ of $G(\mathbb{A}_F^\infty)$ such that $\sigma \otimes \pi_{\infty}$ is isomorphic to a cuspidal automorphic representation of $G(\mathbb{A}_F)$: note that the multiplicities are all 1, which is part of the Jacquet-Langlands theorem, and that when all $n_{\tau}=0$, we must include 1-dimensional automorphic representations that factor through the reduced norm. Since the 1-dimensional constituents do not share Hecke eigenvalues with the infinite-dimensional constituents (by purity), the Jacquet--Langlands correspondence and strong multiplicity one theorem further imply that when we pass to $(M_Q/I_QM_Q) \otimes_{\iota} \C$, the only summand that survives is  
\[
\otimes_{v \not \in S \cup Q} \pi_v^{U_v}/(T_v- \lambda_Q(T_v)) \otimes_{v \in S \setminus \{v \vert \infty\}} \pi_v^{U_v} \otimes_{v \in Q} \pi_v^{U_{Q, v}}/(U_{\vpi_v}-\lambda_Q(U_{\vpi_v})),
\]
%\[
%\otimes_{v \in S \cup Q} \pi_v^{U_{Q, v}}/(U_{\vpi_v}-\lambda_Q(U_{\vpi_v}))
%\]
which we must check is one-dimensional. For $v \not \in S \cup Q$, $\pi_v^{U_v}$ is 1-dimensional with $T_v$ acting by $\lambda_Q(T_v)= \lambda(T_v)$. For $v \in S \setminus \{v \vert \infty\}$, $\pi_v^{U_v}$ is one-dimensional by the choice of $U$, 
%isomorphic to an unramified twist of the Steinberg representation with $U_{\vpi_v}$ acting by the scalar $\lambda_Q(U_{\vpi_v})$ (see \textit{loc. cit.}), 
while for $v \in Q$, $\pi_v^{U_{Q_v}}/(U_{\vpi_v}-\alpha_v)$ is one-dimensional since the $\alpha_v$-eigenspace in $\pi_v^{U_{Q_v}}$ is one-dimensional (whether or not $U_{\vpi_v}$ has distinct eigenvalues).
\end{proof}

We end this section by summarizing the commutative algebra of congruence modules that we will require. We define the congruence module as in \cite[\S 4.4, p.~125, Lemma~4.17]{ddt}, with some broader context taken from \cite{iyengar-khare-manning:congruence-modules}, and the congruence ideal as in
\cite{iyengar-khare-manning:congruence-modules}. Let $\mu\colon T\to\mc{O}$ be an augmentation of a complete Noetherian local
$\mc{O}$-algebra, and put $J=\ker(\mu)$.  Assume that the localization $T_J$ is regular of dimension 0, i.e.,
$(T,\mu)\in\mc{C}_{\mc{O}}(0)$ in the sense of  \cite{iyengar-khare-manning:congruence-modules}. Equivalently, $T_J\simeq E$.  Let $L$ be a finitely generated $T$-module of positive depth such that $L_J\neq 0$. In our work, $L$ is finite flat over $\mc{O}$, which  implies the positive depth and  finite generation hypothesis.

The congruence module of $L$ at $\mu$ is 
\begin{equation}\label{congruence-module}
\Psi_\mu(L):=
\frac{L}{L[J]+L[\Ann_T(J)]}.
\end{equation}
In the self-dual finite-flat setting, this is the module appearing in
\cite[\S~4.4, Lemma~4.17, p.~125]{ddt}; by \cite[Proposition~2.10]{iyengar-khare-manning:congruence-modules}, it is the codimension zero case of the more general construction in \textit{loc. cit.}.

Regarding $\mc{O}$ as a $T$-module through $\mu$, define the congruence ideal
of $L$ by
\begin{equation}\label{module-congruence-ideal}
\eta_\mu(L):=
\im\left(
\Hom_T(\mc{O},L)\otimes_{\mc{O}}\Hom_T(L,\mc{O})
\longrightarrow\mc{O}
\right),
\end{equation}
where the map is given by composition.  This is
\cite[\S~2, equation~(2.4)]{iyengar-khare-manning:congruence-ideals},
specialized to codimension zero.  
For $L=T$, it recovers the usual congruence ideal
\[
\eta_\mu(T)=\mu(\Ann_T(J)).
\]

More precisely, write
\[
\Psi_\mu(L)\simeq
\bigoplus_{i=1}^{d}\mc{O}/\vpi^{e_i}\mc{O},
\qquad
0\leq e_1\leq\cdots\leq e_d,
\]
where $d=\dim_E L_J$.  Then
\[
\eta_\mu(L)=\vpi^{e_1}\mc{O},
\qquad
\Fitt_{\mc{O}}\bigl(\Psi_\mu(L)\bigr)
=\vpi^{e_1+\cdots+e_d}\mc{O};
\]
see \cite[Lemma~2.19]{iyengar-khare-manning:congruence-ideals}.
Thus the congruence ideal records the least elementary divisor of the
congruence module.  In particular, when $d=1$,
\begin{equation}\label{rank-one-congruence}
\Psi_\mu(L)\simeq\mc{O}/\eta_\mu(L).
\end{equation}
We will only need the case $d=1$ and the following consequence: if $L$ carries a perfect $\mc{O}$-linear pairing $\langle \cdot, \cdot \rangle$ for which $T$ acts
self-adjointly, and $\varphi$ spans $L[J]$,  then \eqref{module-congruence-ideal} is the image of the
restriction of this pairing to the $\mu$-eigenspace, hence $\Psi_{\mu}(L) \cong \mc{O}/(\langle \varphi, \varphi \rangle)$.

\section{Taylor-Wiles primes}\label{twsection}
We now make the more precise hypotheses on $\Pi_0$ that we will use for our main analysis. See \S \ref{generalsection} for discussion of how to relax these hypotheses.
\begin{assumption}\label{twassumption}
Assume that $\Pi_0$, $F_0$, and $p$ satisfy the following:
\begin{itemize}
\item $[F_0(\zeta_p):F_0] >2$ (in particular, $p \geq 5$).
\item $\br_{\Pi_0, \iota}|_{\gal{F_0(\zeta_p)}}$ is absolutely irreducible.%if $p=5$ and $\br$ has projective image $\mr{PGL}_2(\mathbb{F}_5)$, then $F(\zeta_5)$ is not contained in $\ov{F}^{\ker(\mr{proj}(\br))}$.
%\item For all $n \geq 1$, there exists $\sigma \in \gal{F, S}$ such that $\sigma|_{F(\zeta_{p^n})}= \id$ and $\bar{\rho}_{\Pi, \iota}(\sigma)$ has distinct eigenvalues.
\item If $p=5$ and $\br_{\Pi_0, \iota}$ has projective image $\mr{PGL}_2(\mathbb{F}_5)$, then $F_0(\zeta_5)$ is not contained in $\ov{F}_0^{\ker(\mr{proj}(\br_{\Pi_0, \iota}))}$.
\end{itemize}
\end{assumption}
Our analysis will require working with automorphic forms on definite quaternion algebras. To descend general $\Pi_0$ via the Jacquet-Langlands correspondence to a definite quaternion algebra, we have to make a base-change first.
\begin{lemma}\label{bc}
There exists a quadratic Galois extension $F/F_0$ such that:
\begin{itemize}
\item $F$ is totally real.
\item Every place in $S_0$ splits in $F$.
\item $[F(\zeta_p):F]>2$.
\item $\br_{\Pi_0, \iota}|_{\gal{F(\zeta_p)}}$ remains absolutely irreducible, and, when $p=5$, still satisfies the last item in Assumption \ref{twassumption}.
\end{itemize}
\end{lemma}
\begin{proof}
This follows from the Grunwald--Wang theorem.  In addition to imposing
the stated local splitting conditions, we choose $F$ linearly disjoint
from the compositum of $F_0(\zeta_p)$ and the field cut out by
$\br_{\Pi_0,\iota}$, using a local condition at an auxiliary prime
outside $S_0$.  This preserves the condition $[F(\zeta_p):F]>2$, the
absolute irreducibility of the residual representation, and the
exceptional condition when $p=5$.
\end{proof}
For the rest of the paper, if $[F_0: \Q]$ is even we take $F=F_0$, and otherwise we now fix an $F$ as in Lemma \ref{bc}. We then make the constructions of \S \ref{generalsetup} and adopt its notation. (For instance, we let $\Pi= \mr{BC}_{F/F_0}(\Pi_0)$, and we descend $\Pi$ to a cuspidal automorphic representation $\pi$ of $G(\mathbb{A}_F)$, where $G=G^D$ is associated to the quaternion algebra $D$ over $F$ ramified exactly at $\{v \vert \infty\}$.) For shorthand we set $\rho= \rho_{\Pi, \iota}= \rho_{\Pi_0, \iota}|_{\gal{F}}$, valued in $\mr{GL}_2(\mc{O})$, and we let $\br= \br_{\Pi, \iota}$ be its reduction.
\begin{lemma}\label{conductor}
For a fixed conductor $f_v$, there are finitely many inertial types arising from Weil-Deligne representations of fixed dimension of conductor at most $f_v$. In particular, for each $v \vert p$, there is a finite extension $F'_v/F_v$ such that all automorphic representations $\sigma$ of $G(\mathbb{A}_F)$ contributing to $M_Q$ for any auxiliary set $Q$ have associated Galois representations $\rho_{\sigma, \iota}$ that become semistable over $F'_v$.
\end{lemma}
\begin{proof}
The first claim follows directly in dimension two.  Since $p\neq2$, the Weil--Deligne representations in question are sums of characters, twists of the Steinberg parameter, or inductions of characters from quadratic extensions of $F_v$; in each case local class field theory gives only finitely many inertial types of bounded conductor. For the second claim, for each $v \vert p$, the conductor subgroup $U_v$ of $\Pi_v$ bounds the conductor of all Weil--Deligne representations at $v$ that can arise from automorphic representations $\sigma$ with level $U_Q$ (a group with local component $U_v$ at $v$). There are finitely many inertial types with this bounded conductor, and they all factor through some finite extensions $F'_v/F_v$. It follows from local--global compatibility at $v \vert p$ that $\rho_{\sigma, \iota}|_{\gal{F_v}}$ becomes semistable over $F'_v$. 
\end{proof}
\begin{rmk}
This lemma will be needed to define the maps $R^{\mr{pst}}_{\br, S \cup Q_n} \to \mathbb{T}_{Q_n}$ below. As an alternative to the following discussion (and to the above lemma), we could have chosen the extension $F/F_0$ to make $\Pi_v$ either unramified or an unramified twist of the Steinberg representation at all $v \vert p$. We would then have $U_v= \mr{GL}_2(\mc{O}_{F_v})$ or $U_0(v)$ for all $v \vert p$, and the potentially semistable (becoming semistable over $F'_v$) deformation ring used below could simply be replaced with a semistable deformation ring. This would work just as well for proving the eventual vanishing of $H^1_g(\gal{F_0, S}, \ad^0(\rho) \otimes E)$, but for the quantitative bound on the torsion analogue we think this setup is more transparent: for instance, when $[F_0: \Q]$ is even, it allows us always to take $F=F_0$.   
\end{rmk}

We next fix integers $a < b$ such that all (labeled) Hodge-Tate weights of $\rho_{\gal{F_v}}$ (for all $v \vert p$) lie in $[a, b]$. We can then consider the category $\mc{E}^{\mr{pst}}_{F, S}$ of finite $\Z_p[\gal{F, S}]$-modules whose restrictions to $\gal{F_v}$ for all $v \vert p$ are isomorphic to quotients of lattices in potentially semistable representations, becoming semistable over $F'_v$, with Hodge-Tate weights in $[a, b]$. Note that this full subcategory of all finite $\Z_p[\gal{F, S}]$-modules is stable under isomorphisms, subquotients, and finite direct sums. The closure properties just noted, together with \cite[Proposition~1.1]{ramakrishna:thesis}, show that the following definition makes sense.
\begin{defn}\label{defring}
Let $R^{\mr{pst}}_{\br,S}$ be the quotient of the universal deformation
ring of
\[
\br\colon\gal{F,S}\longrightarrow\mr{GL}_2(k)
\]
characterized as follows.  For $A\in\mc{C}_{\mc{O}}$, an $A$-valued
deformation $\tau$ factors through $R^{\mr{pst}}_{\br,S}$ if and only if
\begin{itemize}
\item $\det(\tau)=\kappa^{-1}\chi_\iota$;
\item for every $j\geq1$, the underlying finite
$\mathbf Z_p[\gal{F,S}]$-module of
$\tau\bmod\mathfrak m_A^j$ belongs to
$\mc E^{\mr{pst}}_{F,S}$.
\end{itemize}
We caution the reader that $R^{\mr{pst}}_{\br,S}$, like
$\mc E^{\mr{pst}}_{F,S}$, depends on $[a,b]$ and on
$(F'_v)_{v\mid p}$, although these have been omitted from the notation.
\end{defn}

For notational simplicity, in what follows we let $W= \ad^0(\rho_{\Pi, \iota})$, $W_E= W \otimes_{\mc{O}} E$, and $W_{E/\mc{O}}= W \otimes_{\mc{O}} E/\mc{O}$. We also write $W_{m}$ for $\ad^0(\rho_{\Pi, \iota}) \otimes_{\mc{O}} \mc{O}/\vpi^m$, and set $W_m^*=\Hom_{\mc O}(W_m,E/\mc O(1))$. For each integer $m \geq 1$ and finite set $T \supset S$, we consider the following Selmer condition $\mc{D}_{T, m}$:
\begin{itemize}
\item For all $v \in T \setminus \{v \vert p\}$, $\mc{D}_{T, m, v}= H^1(\gal{F_v}, W_{m})$.
\item For $v\mid p$, $\mc{D}_{T,m,v}$ consists of the classes
$c\in H^1(\gal{F_v},W_m)$ for which the associated first-order
deformation
\[
\rho_{m,c}\colon\gal{F_v}\longrightarrow
\mr{GL}_2\bigl((\mc{O}/\vpi^m)[\epsilon]/(\epsilon^2)\bigr)
\]
has determinant $\det(\rho\bmod\vpi^m)$ and whose underlying finite
$\mathbf Z_p[\gal{F_v}]$-module is a quotient of a $\gal{F_v}$-stable lattice in
a potentially semistable representation that becomes semistable over
$F'_v$ and has Hodge--Tate weights in $[a,b]$.
\item For all $v \not \in T$, $\mc{D}_{T, m, v}= H^1_{\mr{unr}}(\gal{F_v}, W_m)$, the unramified classes.
\end{itemize}
Thus the Selmer group $H^1_{\mc{D}_{T, m}}(F, W_m)$ classifies self-extensions $U$ of $\rho \otimes \mc{O}/\vpi^m$, as $\mc{O}/\vpi^m[\gal{F, T}]$-module, that lie in $\mc{E}^{\mr{pst}}_{F, T}$ and have the fixed determinant condition. Note that for all $m \leq m'$, under the reduction modulo $\vpi^m$ map $W_{m'} \to W_m$, $\mc{D}_{T, m', v}$ maps to $\mc{D}_{T, m, v}$ in all cases. Similarly, the maps $W_{m} \xrightarrow{\vpi^{m'-m}} W_{m'}$ induce maps $\mc{D}_{T, m, v} \to \mc{D}_{T, m', v}$, and the preimage of $\mc{D}_{T, m', v}$ under the induced map $H^1(\gal{F_v}, W_m) \to H^1(\gal{F_v}, W_{m'})$ is $\mc{D}_{T, m, v}$. We have similar maps for the dual Selmer conditions $\mc{D}_{T, m, v}^\perp$ (the annihilators of the $\mc{D}_{T, m, v}$ under local duality), but see the proof of Proposition \ref{mainbound} for a subtlety about the preimages.

We are concerned with controlling the Selmer group (compare \cite[Lemma 2.40]{ddt})
\[
H^1_{\mc{D}_S}(F, W_{E/\mc{O}}) = \varinjlim_m H^1_{\mc{D}_{S, m}}(F, W_m) \xleftarrow{\sim} \varinjlim_m \Hom_{\mc{O}}(\mf{p}/\mf{p}^2, \mc{O}/\vpi^m).
\]
Analogous isomorphisms hold with $\mc{D}_{S \cup Q}$ in place of $\mc{D}_S$, and $\mf{p}_Q$ in place of $\mf{p}$, for any Taylor-Wiles set $Q$. To simplify our Selmer group notation, we will often write $H^1_{\mc{D}_S}(F, W_m)$ for $H^1_{\mc{D}_{S, m}}(F, W_m)$. 

For each $n$, let $Q_n$ be a finite set of primes such that for all $v \in Q_n$, $N(v) \equiv 1 \pmod{p^n}$, and $\br(\Frob_v)$ has distinct eigenvalues; we refer to any such set as a set of Taylor--Wiles primes of level $n$.  As in \S \ref{generalsetup} we choose a root $\alpha_v$ of $X^2-\lambda(T_v)X+N(v)\chi(\vpi_v)$, letting $\beta_v$ be the other root, and we let $\mf{m}_{Q_n}=\ker(\lambda_{Q_n}\bmod\vpi)$.  Thus $\mf{m}_{Q_n}$ is the maximal ideal of $\mathbb{T}_{Q_n}(U_{Q_n})$ generated by $\vpi$, $(T_v-\tr(\rho_{\Pi, \iota})(\Frob_v))_{v \not \in S \cup Q_n}$, and $(U_{\vpi_v}-\alpha_v)_{v \in Q_n}$. 

The irreducibility of $\br$ implies that the universal modular pseudodeformation gives a universal deformation 
\[
\rho_{\mf{m}_{Q_n}} \colon \gal{F, S \cup Q_n} \to \mr{GL}_2(\mathbb{T}_{Q_n})
\]
and thus by Lemma \ref{conductor} a homomorphism $R_{\br, S \cup Q_n}^{\mr{pst}} \to \mathbb{T}_{Q_n}$, which by construction surjects onto the reduced subalgebra $\mathbb{T}_{Q_n}^{S \cup Q_n}$ generated by the good Hecke operators $T_v$, $v \not \in S \cup Q_n$. It follows however from the assumption $\alpha_v \not \equiv \beta_v \pmod{\vpi}$ for $v \in Q_n$ that $U_{\vpi_v}$ is in the image of $R_{\br, S \cup Q_n}^{\mr{pst}}$ for all $v \in Q_n$ (see \cite[Corollary 1.8]{taylor:meromorphic}). Thus $R^{\mr{pst}}_{\br, S \cup Q_n} \to \mathbb{T}_{Q_n}$ is surjective. Furthermore, restricting the universal deformation $\rho^{\univ}_{S \cup Q_n}$ to $\gal{F_v}$, $v \in Q_n$, 
\[
\rho^{\univ}_{S \cup Q_n}|_{\gal{F_v}} \sim \begin{pmatrix} \psi_{\alpha_v} & 0 \\ 0 & \psi_{\beta_v} \end{pmatrix}
\]
for a unique choice of characters $\psi_{\alpha_v}$, $\psi_{\beta_v}$ with $\psi_{\alpha_v}(\Frob_v) \equiv \alpha_v \pmod{\vpi}$. In the usual manner we obtain from local class field theory a homomorphism 
\[
\Delta_v \xrightarrow{\psi_{\alpha_v}} (R_{\br, S \cup Q_n}^{\mr{pst}})^\times \to \mathbb{T}_{Q_n}^\times, 
\]
and the fundamental consequence of local-global compatibility at primes in $Q_n$ is that composed with the action of $\mathbb{T}_{Q_n}$ on $M_{Q_n}$, this agrees with the diamond operator action of $\Delta_v$ on $M_{Q_n}$ (\cite[Corollary 1.8]{taylor:meromorphic}).
 
We summarize our modular input for Bloch-Kato vanishing:
\begin{prop}\label{twsetup}
\begin{enumerate} For any set $Q_n$ as above:
\item The $\mc{O}$-module $I_{Q_n}/I_{Q_n}^2$ is finitely-generated with an upper bound $g$ on the number of generators depending only on $|Q_n|$.
\item $M_{Q_n}$ is a free $\mc{O}[\Delta_{Q_n}]$-module of rank equal to the $\mc{O}$-rank of $S(U, \mc{O})_{\mf{m}}$, and $M_{Q_n}/\mf{a}_{Q_n}M_{Q_n} \xrightarrow{\sim} S(U, \mc{O})_{\mf{m}}$.
\item The augmentation ideal $\mf{a}_{Q_n}$ is contained in $I_{Q_n}$. The $\mc{O}$-module $M_{Q_n}/I_{Q_n}M_{Q_n}$ is independent of $Q_n$, isomorphic to $M_{\emptyset}/IM_{\emptyset}= S(U, \mc{O})_{\mf{m}}/I S(U, \mc{O})_{\mf{m}}$. Thus it is isomorphic to $\mc{O} \oplus A_{\emptyset}$ for some finite torsion $\mc{O}$-module $A_{\emptyset}$, %where $A_{Q_n}$ is a finite torsion module whose length is bounded independently of $Q_n$ by 
and we set $a_{\emptyset}:= \lg(A_{\emptyset})= \lg\bigl((S(U, \mc{O})_{\mf{m}}/I S(U, \mc{O})_{\mf{m}})_{\mr{tor}}\bigr)$.
\item The induced map
\[
\mf{a}_{Q_n}M_{Q_n}/\mf{a}_{Q_n}^2 M_{Q_n} \to I_{Q_n}M_{Q_n}/I_{Q_n}^2 M_{Q_n}
\]
has cokernel with length %bounded by a constant 
$b_{\emptyset}= \lg(IM_{\emptyset}/I^2M_{\emptyset})$ independent of $Q_n$.
\item The congruence module $\Psi_{\lambda_{Q_n}}(M_{Q_n})$ of
\eqref{congruence-module} is cyclic and satisfies
\[
\lg\bigl(\Psi_{\lambda_{Q_n}}(M_{Q_n})\bigr)
=
\lg\bigl(\Psi_{\lambda}(M_{\emptyset})\bigr)
+
\lg(\mc{O}/|\Delta_{Q_n}|).
\]
In particular,
\[
\lg\bigl(\Psi_{\lambda_{Q_n}}(M_{Q_n})\bigr)
\geq
\lg(\mc{O}/|\Delta_{Q_n}|)
=e\,v_p(|\Delta_{Q_n}|)\geq ne|Q_n|.
\]
\end{enumerate}
\end{prop}
\begin{proof}
(1): Let $\mf{p}_{Q_n}$ be the kernel of the composite $R^{\mr{pst}}_{\br, S \cup Q_n} \to \mathbb{T}_{Q_n} \xrightarrow{\lambda_{Q_n}} \mc{O}$. There exists $g \in \Z$ depending only on $\br$ and $|Q_n|$ such that there is a surjection $\mc{O}[[X_1, \ldots, X_g]] \onto R^{\mr{pst}}_{\br, S \cup Q_n}$: indeed, we can take 
\[
g= \dim_k H^1_{\mc{D}_{S \cup Q_n}}(F, \ad^0(\br)) \leq \dim_k H^1_{\mc{D}_S}(F, \ad^0(\br))+ |Q_n|,
\]
since $Q_n$ being a set of Taylor-Wiles primes ensures $\dim H^1(\gal{F_v}, \ad^0(\br))/ H^1(\gal{F_v}/I_{F_v}, \ad^0(\br)) =1$ for all $v \in Q_n$.
The $\mc{O}$-module surjection $\mf{p}_{Q_n}/\mf{p}_{Q_n}^2 \onto I_{Q_n}/I_{Q_n}^2$ then implies $I_{Q_n}/I_{Q_n}^2$ is generated by at most $g$ elements.

(2): By our assumption $[F(\zeta_p):F]>2$, the level subgroup $U$ is sufficiently small in the sense that $p$ does not divide the order of $(U (\mathbb{A}_F^{\infty})^\times \cap g^{-1} D^\times g)/F^\times$ for any $g \in (D \otimes \mathbb{A}_F^{\infty})^\times$: see the proof of \cite[Lemma 1.1]{taylor:meromorphic}. Part (2) is then a standard consequence of the $Q_n$ being sets of Taylor-Wiles primes: see \cite[Lemma 2.3]{taylor:meromorphic}.

(3): We have noted above that the diamond operators $\mc{O}[\Delta_v]$ are in $\mathbb{T}_{Q_n}$, and $\delta \in \Delta_v$ acts trivially on the $U_{\vpi_v}$-eigenform (with eigenvalue $\alpha_v$) inducing $\lambda_{Q_n}$ (see the calculation leading to \cite[Fact 5.7]{gee:aws}). % by $\chi_{\alpha_v}(\delta)$, where $\Pi_v \cong \ind; 
Thus $\mf{a}_{Q_n}$ is contained in $I_{Q_n}$. Moreover,
\[
M_{Q_n}/I_{Q_n}M_{Q_n}
\xrightarrow{\sim}
(M_{Q_n}/\mf{a}_{Q_n}M_{Q_n})/(I_{Q_n}/\mf{a}_{Q_n})(M_{Q_n}/\mf{a}_{Q_n}M_{Q_n})
\xrightarrow{\sim}
S(U, \mc{O})_{\mf{m}}/I S(U, \mc{O})_{\mf{m}},
\]
where $I=\ker(\lambda)$. In particular, this last quotient is a finite $\mc{O}$-module of generic rank 1 (Lemma \ref{weakgenfree}) that is independent of $Q_n$.

(4): The cokernel in question is isomorphic to $I_{Q_n}(M_{Q_n}/\mf{a}_{Q_n})/I_{Q_n}^2(M_{Q_n}/\mf{a}_{Q_n})$; since $Q_n$ is a Taylor-Wiles set of primes, as noted above $\mathbb{T}_{Q_n}$ equals the reduced subalgebra $\mathbb{T}^{S \cup Q_n}_{Q_n}$, and, by part (2), the cokernel is isomorphic, as a $\mathbb{T}^{S \cup Q_n}_{Q_n} \onto \mathbb{T}_{\emptyset}$-module, to $I M_{\emptyset}/I^2 M_{\emptyset}$. As such, it is independent of $Q_n$, and it is finite since it is a finitely generated module over $\mathbb{T}_{\emptyset}/(I+ \Ann_{\mathbb{T}_{\emptyset}}(I))$, which has finite cardinality since $\mathbb{T}_{\emptyset}$ is reduced.

(5): We use the congruence formalism of \eqref{congruence-module}--\eqref{rank-one-congruence}. The three modules $M_Q$, $M_Q^{\mr{St}}:= S(U_0(Q), \mc{O})_{\mf{m}_Q}$, and $M_{\emptyset}$ carry perfect $\mc{O}$-linear pairings $\langle \cdot , \cdot \rangle_{V}$, for $V= U_Q, U_0(Q), U$, for which the Hecke algebras $\mathbb{T}_Q$, $\mathbb{T}_Q^{\mr{St}}$, and $\mathbb{T}$ consist of self-adjoint operators (see \cite[p. 742, 752]{taylor:meromorphic}). More precisely, recall from \S \ref{generalsetup} the perfect pairing $\langle \cdot, \cdot \rangle$ on the codomain of our integral automorphic forms. We obtain perfect pairings $\langle \cdot, \cdot \rangle_{0, V}$ on each $S(V, \mc{O})$ using the formula from \cite[p.742]{taylor:meromorphic}, and \textit{loc. cit.} also observes that the adjoint of the Hecke operator $[VgV]$ with respect to this pairing is $\chi_{\iota}(g)[V g^{-1}V]$. It follows easily that for $v \not \in S \cup Q$, $T_v$ is self-adjoint with respect to $\langle \cdot, \cdot \rangle_{0, V}$. Since we also want this self-adjointness for the operators $U_{\vpi_v}$, $v \in Q$, we instead use the pairing
\[
\langle \varphi_1, \varphi_2 \rangle_V= \prod_{v \in Q} \chi_{\iota}(\vpi_v)^{-1} \langle \varphi_1, w_Q \varphi_2 \rangle_{0, V},
\]
where $w_Q \in \mr{GL}_2(\mathbb{A}_F^\infty)$ has component $\begin{pmatrix} 0 & 1 \\ \vpi_v & 0 \end{pmatrix}$ for each $v \in Q$, and component 1 elsewhere. For this modified pairing, the full Hecke algebras $\mathbb{T}_Q$ and $\mathbb{T}_Q^{\mr{St}}$ are self-adjoint, and in particular we get induced perfect pairings on $M_Q$, $M_Q^{\mr{St}}$, and $M_{\emptyset}$ to which we can apply the formalism of \eqref{congruence-module}--\eqref{rank-one-congruence}.

We have by part (2) (as explained in more detail in \cite[Lemmas~2.2--2.3 and Corollary~2.4]{taylor:meromorphic}), $\mathbb{T}_Q \onto \mathbb{T}_Q^{\mr{St}} \onto \mathbb{T}$-equivariant maps $\beta= \beta_2 \circ \beta_1$
\[
\beta \colon M_Q \xrightarrow{\beta_1} M_Q^{\mr{St}} \xrightarrow[\sim]{\beta_2} M_{\emptyset}
\]
where $\beta_1$ is the trace, and $\beta_2^{-1}(f)= \prod_{v \in Q} \left( A_v - \begin{pmatrix} 1 & 0 \\ 0 & \vpi_v \end{pmatrix}  \right) f$, where $A_v \in \mathbb{T}$ is the unique root of $X^2-T_vX+ N(v)\chi_{\iota}(\vpi_v)$ lifting $\alpha_v$. For any of these maps $\beta_{?}$, denote by $\beta_{?}^{\dagger}$ the adjoint with respect to the given pairings. Let $\varphi$ span $M_{\emptyset}[I]$, and let $\varphi_Q$ span $M_Q[I_Q]$. Since $\beta$ is surjective, $\coker(\beta^\dagger)$ is torsion-free, and thus $\beta^\dagger(\varphi)$ equals $\varphi_Q$ up to a unit $u \in \mc{O}^\times$. We now compute 
\[
\langle \varphi_Q, \varphi_Q \rangle_{U_Q}= u^{-2} \langle \beta^\dagger \varphi, \beta^{\dagger} \varphi \rangle_{U_Q}= u^{-2}\langle \varphi, \beta \circ \beta^{\dagger} \varphi \rangle_U.
\]
We can compute $\beta \beta^\dagger= \beta_2 \beta_1 \beta_1^{\dagger} \beta_2^\dagger= |\Delta_Q| \beta_2 \beta_2^\dagger$, since the definition of the pairings implies $\beta_1^{\dagger}$ is just the inclusion, and each diamond operator acts trivially on $\beta^{\dagger}(\varphi)$. Finally, $\beta_2$ is an isomorphism, so $\beta_2^\dagger$ is as well, and 
%$\beta_2 \beta_2^\dagger$ is multiplication by a unit, by the calculation in \cite[Lemma 2]{taylor:hmf} and the fact that $Q$ consists of Taylor-Wiles primes. 
we conclude that $\lg(\Psi_{\lambda_Q}(M_Q))= \lg(\Psi_{\lambda}(M_{\emptyset}))+ \lg(\mc{O}/|\Delta_Q|)$, which implies (5).

%We also record the output of this argument when $Q$ consists of primes $v$ such that $N(v) \equiv 1 \pmod{p^n}$ but the roots $\alpha_v$ and $\beta_v$ of $X^2-\lambda(T_v)X+N(v)\chi(\vpi_v)$ satisfy $\alpha_v \equiv \beta_v \pmod{\vpi^M}$, $\alpha_v \not \equiv \beta_v \pmod{\vpi^{M+1}}$.
\end{proof}
%In the above proof we used the surjection $\mathbb{T}_{Q_n} \onto \mathbb{T}_{\emptyset}$, which induces a surjection $\alpha \colon I_{Q_n}/I_{Q_n}^2 \onto I/I^2$ with the property that $\ker(\alpha)$ contains the image of $\mf{a}_{Q_n}/\mf{a}_{Q_n}^2$.
\begin{rmk}
In the toy case where $M_{\emptyset}$ is free of rank 1 over $\mathbb{T}_{\emptyset}$, $a_{\emptyset}=0$, and $b_{\emptyset}= \lg(I/I^2)$ is the cotangent space in Wiles's numerical criterion.
\end{rmk}

Here and below, $\mathfrak aM/I^2M$ denotes the image of
$\mathfrak aM$ in $M/I^2M$. We use the following simple  lemma in the proof of
Proposition~\ref{growth}.

\begin{lemma}\label{cokernel-beta}
Let $\lambda\colon T\to\mc O$ be an augmentation of a complete
Noetherian local $\mc O$-algebra, let $I=\ker(\lambda)$, and let
$\mathfrak a\subseteq I$ be an ideal of $T$.  Let $M$ be a finite
$T$-module, and suppose that
\[
M/IM\simeq\mc O x\oplus A,
\qquad A=(M/IM)_{\mr{tor}}.
\]
Then the cokernel of the multiplication map
\[
\mathfrak a/\mathfrak a^2\longrightarrow
\Hom_{\mc O}(M/IM,\mathfrak aM/I^2M)
\]
is annihilated by $\Ann_{\mc O}(A)^2$.
\end{lemma}

\begin{proof}
Let $t$ generate $\Ann_{\mc O}(A)$, and let
\[
f\in\Hom_{\mc O}(M/IM,\mathfrak aM/I^2M).
\]
Then $tf$ vanishes on $A$.  Write
\[
tf(x)=\sum_i b_i m_i\pmod{I^2M},
\qquad b_i\in\mathfrak a,
\]
and write the image of $m_i$ in $M/IM$ as
\[
m_i=c_i x+a_i,\qquad c_i\in\mc O,\quad a_i\in A.
\]
Since $tA=0$, we have
\[
t^2f(x)=\left(\sum_i tc_i b_i\right)x\pmod{I^2M}.
\]
Multiplication by $\sum_i tc_i b_i\in t\mathfrak a$ kills
$A\subset M/IM$ as a map to $\mathfrak aM/I^2M$.  Hence $t^2f$ is
induced by multiplication by $\sum_i tc_i b_i$, and therefore lies
in the image of $\mathfrak a/\mathfrak a^2\longrightarrow
\Hom_{\mc O}(M/IM,\mathfrak aM/I^2M)$.
\end{proof}

\begin{prop}\label{growth}
Fix an integer $r \geq 0$. There is a constant $d_r'$ such that for all $n \geq 1$, and all Taylor--Wiles sets $Q_n$ of level $n$ and size $|Q_n|=r$, the image of $\mf{a}_{Q_n}/\mf{a}_{Q_n}^2 \to I_{Q_n}/I_{Q_n}^2$ has length at least $n e r-d_r'$. In terms of the constants of Proposition \ref{twsetup}, we can take $d_r'= b_{\emptyset}+r a_{\emptyset}$.  

For the constant $d_r$ defined in Theorem \ref{introthm}, the image of
\[
\Hom_{\mc{O}}(I_{Q_n}/I_{Q_n}^2, \mc{O}/\vpi^{ne}) \to \Hom_{\mc{O}}(\mf{a}_{Q_n}/\mf{a}_{Q_n}^2, \mc{O}/\vpi^{ne})
\]
has length at least $\lg (\Hom_{\mc{O}}(\mf{a}_{Q_n}/\mf{a}_{Q_n}^2, \mc{O}/\vpi^{ne}))-d_r= ner-d_r$.% for some bounded error term $b$ depending only on $|Q_n|$.
\end{prop}
\begin{proof}
Using Proposition \ref{twsetup}, we write $M_{Q_n}/I_{Q_n} M_{Q_n}= \mc{O}x \oplus A_{Q_n}$, where $x$ is not torsion, and %$|A_{Q_n}|$ is bounded depending only on $|Q_n|$. 
$A_{Q_n} \cong A_{\emptyset}$ is up to isomorphism independent of $Q_n$. We control the image of $\mf{a}_{Q_n}/\mf{a}_{Q_n}^2 \to I_{Q_n}/I_{Q_n}^2$ by considering their respective images in $\End(M_{Q_n}/I_{Q_n}^2 M_{Q_n})$, 
%which we analyze as submodules of
%\[
%\Hom_{\mc{O}}(M_{Q_n}/I_{Q_n}M_{Q_n}, I_{Q_n}M_{Q_n}/I_{Q_n}^2M_{Q_n}) \xrightarrow{\sim} \Hom_{\mc{O}}(\mc{O}x, I_{Q_n}M_{Q_n}/I_{Q_n}^2M_{Q_n}) \oplus \Hom_{\mc{O}}(A_{Q_n}, I_{Q_n}M_{Q_n}/I_{Q_n}^2M_{Q_n}).
%\]
 via the commutative diagram
\begin{equation}\label{keydiagram}
\xymatrix{
\mf{a}_{Q_n}/\mf{a}_{Q_n}^2 \ar[d]^{\alpha} \ar[r]^-{\beta} & \Hom(M_{Q_n}/I_{Q_n}M_{Q_n}, \mf{a}_{Q_n} M_{Q_n}/ I_{Q_n}^2 M_{Q_n}) \ar[d]^\gamma \\
I_{Q_n}/I_{Q_n}^2 \ar[r]^-{\delta} & \Hom(M_{Q_n}/I_{Q_n} M_{Q_n}, I_{Q_n} M_{Q_n}/I_{Q_n}^2 M_{Q_n}),
}
\end{equation}
where the maps $\beta$ and $\delta$ are induced by the action of $\mathbb{T}_{Q_n}$ on $M_{Q_n}$.

%The image of $\mf{a}_{Q_n}$ in $\End(M_{Q_n}/I_{Q_n}^2 M_{Q_n})$ 
The image of $\beta$ has length at least that of its quotient (by restricting the endomorphism) $\mf{a}_{Q_n}\cdot x \subset \mf{a}_{Q_n} M_{Q_n}/I_{Q_n}^2 M_{Q_n}$. Since $M_{Q_n}/I_{Q_n}M_{Q_n}= \mc{O}x \oplus A_{Q_n}$, 
%, and again using that $M_{Q_n}/I_{Q_n}M_{Q_n}= \mc{O}x \oplus A_{Q_n}$ 
$\lg (\mf{a}_{Q_n} M_{Q_n}/I^2_{Q_n}M_{Q_n}) \leq \lg(\mf{a}_{Q_n} \cdot x)+\lg(\mf{a}_{Q_n} \cdot A_{Q_n})$, so
%The image of $\beta$ has length at least that of its quotient (by restricting the endomorphism) $\mf{a}_{Q_n}\cdot x \subset \mf{a}_{Q_n} M_{Q_n}/I_{Q_n}^2 M_{Q_n}$, and again using that $M_{Q_n}/I_{Q_n}M_{Q_n}= \mc{O}x \oplus A_{Q_n}$ we find that 
\[
\lg(\im(\beta)) \geq \lg (\mf{a}_{Q_n} \cdot x) \geq \lg(\mf{a}_{Q_n}M_{Q_n}/I_{Q_n}^2 M_{Q_n})-c_{|Q_n|},
\]
where $c_{|Q_n|} \leq |Q_n| \lg(A_{Q_n})= |Q_n| a_{\emptyset}$ is bounded by the number of $\mc{O}$-module generators of $\mf{a}_{Q_n}/\mf{a}_{Q_n}^2$ (namely $|Q_n|$) times the length of each such generator's image on $A_{Q_n}$.

%which equals $\lg I_{Q_n} M_{Q_n}/I_{Q_n}^2 M_{Q_n}$ up to an error term $\lg (I_{Q_n}/I_{Q_n}^2) \cdot A_{Q_n}$ that is at most $g \lg(A_{Q_n})$ (bounded depending only on $|Q_n|$; here $g$ is as in Proposition \ref{twsetup}). The same reasoning as above shows that $\lg(\im(\beta))$ is at least $\lg(\mf{a}_{Q_n} M_{Q_n}/I_{Q_n}^2 M_{Q_n})$ up to a constant $a_1$ (that is $O(|{Q_n}|)$). 

Since $\gamma$ is injective, and the cokernel of $\mf{a}_{Q_n}M_{Q_n}/I_{Q_n}^2 M_{Q_n} \to I_{Q_n}M_{Q_n}/I_{Q_n}^2 M_{Q_n}$ has length bounded by $b_{\emptyset}$ (by part (4) of Proposition \ref{twsetup}), we conclude that  $\lg(\im(\gamma \circ \beta)) \geq \lg(I_{Q_n}M_{Q_n}/I_{Q_n}^2 M_{Q_n})-b_{\emptyset}-c_{|Q_n|}$. 
%\geq ne |{Q_n}|+a_1+a_2$. 
We can bound this in terms of the congruence ideal using basic properties of Fitting ideals:
\begin{align*}
\Fitt_{\mc{O}}(I_{Q_n} M_{Q_n}/I_{Q_n}^2 M_{Q_n})&= \Fitt_{\mathbb{T}_{Q_n}/I_{Q_n}}(I_{Q_n} M_{Q_n} \otimes_{\mathbb{T}_{Q_n}} \mathbb{T}_{Q_n}/ I_{Q_n})= \Fitt_{\mathbb{T}_{Q_n}}(I_{Q_n} M_{Q_n}) \otimes_{\mathbb{T}_{Q_n}} \mathbb{T}_{Q_n}/I_{Q_n} \\ &\subset \Ann_{\mathbb{T}_{Q_n}}(I_{Q_n} M_{Q_n}) \otimes_{\mathbb{T}_{Q_n}} \mathbb{T}_{Q_n}/I_{Q_n}= \lambda_{Q_n}(\Ann_{\mathbb{T}_{Q_n}}(I_{Q_n} M_{Q_n})) \\
&\subseteq \Ann_{\mc{O}}\bigl(\Psi_{\lambda_{Q_n}}(M_{Q_n})\bigr)
=\eta_{\lambda_{Q_n}}(M_{Q_n}).
\end{align*}
Thus by part (5) of Proposition \ref{twsetup}, we conclude (recall from part (5) of Proposition \ref{twsetup} that $\Psi_{\lambda_{Q_n}}(M_{Q_n})$ is cyclic) that
\[
\lg(\im(\gamma \circ \beta)) \geq \lg(\Psi_{\lambda}(M_{\emptyset}))+ \lg(\mc{O}/|\Delta_{Q_n}|)-b_{\emptyset}- |Q_n|a_{\emptyset},
\]
and since $\gamma \circ \beta$ factors through $\alpha$ we obtain the same lower bound for $\lg(\im(\alpha))$. This proves the claim about $d_r'$. Since the source of $\gamma \circ \beta$ has length $\lg(\mc{O}/|\Delta_{Q_n}|)$, we also find that $\lg \ker(\gamma \circ \beta) \leq b_{\emptyset}+|Q_n|a_{\emptyset}-\lg(\Psi_{\lambda}(M_{\emptyset}))$.%, a bound independent of $Q_n$.

For the second claim, we take the Pontryagin dual (denoted $(?)^*$) of the diagram \eqref{keydiagram}. By exactness of the Pontryagin dual, $\lg(\im(\beta^* \circ \gamma^*))= \lg(\im(\gamma \circ \beta))$. Moreover, $\coker(\gamma)$ is annihilated by $\Ann(IM_{\emptyset}/I^2M_{\emptyset})$ by part (4) of Proposition \ref{twsetup}, while Lemma \ref{cokernel-beta} shows that $\coker(\beta)$ is annihilated by $\Ann(A_{\emptyset})^2$. Thus $\coker(\gamma \circ \beta)$ is killed by
\[
\vpi^{j_{\emptyset}}\mc{O}
=\Ann(IM_{\emptyset}/I^2M_{\emptyset})\,
 \Ann(A_{\emptyset})^2.
\]

We thus have an exact sequence
\[
0 \to K_{Q_n} \to \mf{a}_{Q_n}/\mf{a}_{Q_n}^2 \xrightarrow{\gamma \circ \beta} \Hom(M_{Q_n}/I_{Q_n}M_{Q_n}, I_{Q_n}M_{Q_n}/I^2_{Q_n}M_{Q_n}) \to L_{Q_n} \to 0
\]
where $\lg(K_{Q_n})$ is bounded in a way depending only on $|Q_n|$, and the exponent of $L_{Q_n}$ is bounded independently of $Q_n$ (by $j_{\emptyset}$). If $ne\leq j_{\emptyset}$, the asserted lower bound is automatic.
We may therefore assume that $ne>j_{\emptyset}$. Now, any element of 
\[
\ker((\mf{a}_{Q_n}/\mf{a}_{Q_n}^2)^* \to K_{Q_n}^*)[\vpi^{ne-j_{\emptyset}}]
\]
is in the image of $(\gamma \circ \beta)^*[\vpi^{ne}]$, the restriction of
$(\gamma \circ \beta)^*$ to the $\vpi^{ne}$-torsion.  Hence
\[
\lg\!\left(\im\bigl((\gamma \circ \beta)^*[\vpi^{ne}]\bigr)\right)
\geq
\lg\!\left(\Hom(\mf{a}_{Q_n}/\mf{a}_{Q_n}^2,
\mc{O}/\vpi^{ne-j_{\emptyset}})\right)-\lg(K_{Q_n}).
\]
In particular, since $\gamma \circ \beta$ factors through $\alpha$, the image of
\[
\Hom(I_{Q_n}/I_{Q_n}^2, \mc{O}/\vpi^{ne}) \to \Hom(\mf{a}_{Q_n}/\mf{a}_{Q_n}^2, \mc{O}/\vpi^{ne})
\]
has length at least
\[
(ne-j_{\emptyset})|Q_n|-b_{\emptyset}
-a_{\emptyset}|Q_n|
+\lg\bigl(\Psi_{\lambda}(M_{\emptyset})\bigr)
=ner-d_r,
\]
where $r=|Q_n|$ and
\[
d_r=r(j_{\emptyset}+a_{\emptyset})+b_{\emptyset}
-\lg\bigl(\Psi_{\lambda}(M_{\emptyset})\bigr).
\]

\end{proof}

We now use a variant of Wiles's argument from \cite[Theorem 3.1(i)]{wiles:fermat} to convert Proposition \ref{growth} for Taylor-Wiles sets into a bound on the Selmer group $H^1_{\mc{D}_S}(F, W_{E/\mc{O}}).$ The following central proposition translates Proposition \ref{growth} into the growth of the Selmer groups $H^1_{\mc{D}_{S \cup Q_n}}(F, W_{ne})$.
\begin{cor}\label{galcohgrowth}
For any $r \geq 0$ and any Taylor--Wiles set $Q_n$ of level $n$ and size $|Q_n|=r$,
\[
h^1_{\mc{D}_{S \cup Q_n}}(\gal{F, S \cup Q_n}, W_{ne})- h^1_{\mc{D}_S}(\gal{F, S}, W_{ne}) \geq ner-d_r,
\]
where $d_r$ is the constant in Proposition \ref{growth}.
\end{cor}
\begin{proof}
The $\mc{O}$-algebra homomorphisms $\mc{O}[\Delta_{Q_n}] \xrightarrow{\theta} R^{\mr{pst}}_{\br,S\cup Q_n} \to \mathbb{T}_{Q_n}$ give rise to a commutative diagram
\[
\xymatrix{
\Hom_{\mc{O}}(\mf{p}/\mf{p}^2, \mc{O}/\vpi^{ne}) \ar[r] & \Hom_{\mc{O}}(\mf{p}_{Q_n}/\mf{p}_{Q_n}^2, \mc{O}/\vpi^{ne}) \ar[r]^{\theta^*} & \Hom_{\mc{O}}(\mf{a}_{Q_n}/\mf{a}_{Q_n}^2, \mc{O}/\vpi^{ne}) \\
\Hom_{\mc{O}}(I/I^2, \mc{O}/\vpi^{ne}) \ar[u] & \Hom_{\mc{O}}(I_{Q_n}/I_{Q_n}^2, \mc{O}/\vpi^{ne}) \ar[u] \ar[ur]^{\alpha^*}, &
}
\]
and all we need to know is that the top row is a complex (it is in fact exact) whose first map is injective.. Proposition \ref{growth} shows that $\alpha^*$, and thus $\theta^*$, has image with length at least $ner-d_r$. Reinterpreting in Galois cohomology, the isomorphism
\[
\Hom_{\mc{O}}(\mf{p}_{Q_n}/\mf{p}_{Q_n}^2, \mc{O}/\vpi^{ne}) \xrightarrow{\sim} H^1_{\mc{D}_{S \cup Q_n}}(\gal{F, S \cup Q_n}, W_{ne}),
\]
and analogously when $Q_n= \emptyset$, shows that
\[
h^1_{\mc{D}_{S \cup Q_n}}(\gal{F, S \cup Q_n}, W_{ne})- h^1_{\mc{D}_S}(\gal{F, S}, W_{ne}) \geq ner-d_r.
\]
\end{proof}

\begin{lemma}\label{gw}
There is a constant $c_S$, depending only on the restrictions
$\rho_{\Pi,\iota}|_{\gal{F_v}}$ for $v\in S$, such that, for every
$n\geq1$ and every Taylor--Wiles set $Q_n$ of level $n$,
\[
h^1_{\mc D_{S\cup Q_n}}(F,W_{ne})
-h^1_{\mc D_{S\cup Q_n}^{\perp}}(F,W_{ne}^*)
-ne|Q_n|=c_{S,n},
\]
where $0\leq c_{S,n}\leq c_S$.
\end{lemma}

\begin{proof}
The Greenberg--Wiles formula gives
\begin{align}
&h^1_{\mc D_{S\cup Q_n}}(F,W_{ne})
-h^1_{\mc D_{S\cup Q_n}^{\perp}}(F,W_{ne}^*) \notag\\
&\quad =
h^0(F,W_{ne})-h^0(F,W_{ne}^*)
+\sum_{v\in S\cup Q_n}
\left(
\lg\mc D_{S\cup Q_n,ne,v}
-h^0(\gal{F_v},W_{ne})
\right)
\notag\\
&\quad =
\sum_{v\in S \setminus \{w \vert \infty\}}
\left(
\lg\mc D_{S\cup Q_n,ne,v}
-h^0(\gal{F_v},W_{ne})
\right)
+ne|Q_n|-ne[F:\Q].
\label{gweq}
\end{align}
Here we have used the absolute irreducibility of
$\br|_{\gal{F(\zeta_p)}}$, the oddness of $\br$, and the usual local
duality calculation at the Taylor--Wiles primes.

For each finite place $v\in S$, let $\overline R_v^\square$ denote the following
fixed-determinant lifting ring of $\br|_{\gal{F_v}}$, with determinant
$\kappa^{-1}\chi_\iota$.  If $v\nmid p$, set
\[
\overline R_v^\square=R_{\br|_{\gal{F_v}}}^\square.
\]
If $v\mid p$, let
\[
\overline R_v^\square
=R_{\br|_{\gal{F_v}}}^{\square,\mr{pst}},
\]
the lifting ring whose Artinian points are lifts occurring as
subquotients of lattices in potentially semistable representations
with Hodge--Tate weights in $[a,b]$ that become semistable over $F'_v$. A theorem of Tong Liu \cite{liu:torsion} implies that the generic fiber $\overline{R}_v^{\square}[1/p]$ is isomorphic to the generic fiber of Kisin's potentially semistable deformation ring, which has the universal property of being the quotient of (the fixed determinant lifting ring) $R^{\square}_{\br|_{\gal{F_v}}}$ whose points in finite $E$-algebras correspond to representations that become semistable over $F_v'$ and have all Hodge-Tate weights in $[a, b]$ (see \cite[Theorem 2.7.6]{kisin:pst}).

Let $\mf p_v$ be the kernel of the homomorphism
$\overline R_v^\square\to\mc O$ defined by
$\rho|_{\gal{F_v}}$.  Then
\[
\lg\mc D_{S\cup Q_n,ne,v}
=
\lg\Hom_{\mc O}(\mf p_v/\mf p_v^2,\mc O/\vpi^{ne})
-\bigl(\lg W_{ne}-h^0(\gal{F_v},W_{ne})\bigr).
\]
Moreover,
\[
\Hom_{\mc O}(\mf p_v/\mf p_v^2,\mc O/\vpi^{ne})
\cong
\Hom_{\mc O}(C_v,\mc O/\vpi^{ne}),
\qquad
C_v:=
\Omega_{\overline R_v^\square/\mc O}
\otimes_{\overline R_v^\square,\rho}\mc O.
\]

By local--global compatibility and the genericity of $\Pi_v$, the
Weil--Deligne representation associated to
$\rho|_{\gal{F_v}}$ is generic; see
\cite[Lemma~1.1.3]{allen:selmer}.  Suppose that $v\mid p$, and let
$\tau_v$ and $\mathbf v_v$ denote the Galois type and Hodge type of
$\rho|_{\gal{F_v}}$.  By
\cite[Theorem~1.2.7]{allen:selmer}, the point defined by
$\rho|_{\gal{F_v}}$ is smooth on the generic fiber of the corresponding
fixed-type deformation ring.  It therefore lies on a unique component
of that space.  Since the Hodge and Galois types are uniquely determined
at this point, this is also the unique component of
$\Spec\overline R_v^\square[1/p]$ passing through it.  Thus the completed
local rings agree, and \cite[Lemma~1.2.5]{allen:selmer} identifies the
image of the resulting Zariski tangent space in $H^1(F_v,W_E)$ with
$H^1_g(F_v,W_E)$.  For $v\nmid p$, the corresponding assertion follows
from \cite[Proposition~1.2.2]{allen:selmer}.  Consequently, $C_v$ has rank
\[
d_v=
\begin{cases}
3,&v\nmid p,\\
3+[F_v:\Q_p],&v\mid p.
\end{cases}
\]
Write
\begin{equation}\label{localtorsion}
C_v\simeq
\mc O^{d_v}\oplus\bigoplus^{t_v}_{i=1}\mc O/\vpi^{e_{v,i}}.
\end{equation}
It follows that
\[
\lg\mc D_{S\cup Q_n,ne,v}
-h^0(\gal{F_v},W_{ne})
=
\begin{cases}
c_{v,n},&v\nmid p,\\
c_{v,n}+ne[F_v:\Q_p],&v\mid p,
\end{cases}
\]
where
\[
c_{v,n}=\sum_i\min\{ne,e_{v,i}\},
\qquad
0\leq c_{v,n}\leq c_v:=\sum_i e_{v,i}.
\]
Since $S$ contains all places above $p$,
\[
\sum_{v\mid p}[F_v:\Q_p]=[F:\Q].
\]
Equation~\eqref{gweq} therefore gives
\[
c_{S,n}=\sum_{v\in S}c_{v,n},
\qquad
c_S=\sum_{v\in S}c_v,
\]
which proves the lemma.
\end{proof}

Now by judicious choice of the sets $Q_n$, we obtain the main bound:
\begin{prop}\label{mainbound}
There is an integer $r\geq0$ such that, for every $n\geq1$, there is a
set $Q_n$ of Taylor--Wiles primes of level $n$, with $|Q_n|=r$, such that
\[
h^1_{\mc D_S}(F,W_{ne})
\leq d_r+c_{S,n}
\leq d_r+c_S,
\]
where $d_r$ is the constant appearing in Proposition~\ref{growth},
and $c_{S,n}$ and $c_S$ are the local error terms defined in
Lemma~\ref{gw}.
\end{prop}
\begin{proof}
We give two explanations, the first following \cite[Theorem 3.1]{wiles:fermat} and the second giving the stated refinement.
Since $\br|_{\gal{F(\zeta_p)}}$ is absolutely irreducible, $(\ad^0(\br)^*)^{\gal{F}}$ vanishes, and inductively $(W_m^*)^{\gal{F}}=0$ for all $m \geq 1$. Thus for any $n \geq 1$, the map (induced from the kernel of multiplication by $\vpi$) $H^1(\gal{F, S \cup Q}, W_1^*) \to H^1(\gal{F, S \cup Q}, W_{ne}^*)[\vpi]$ is an isomorphism (for any $Q$). 
%Moreover, the preimage of $\mc{D}_{S , ne, v}^{\perp}$ under $W_1^* \to W_{ne}^*$ is $\mc{D}_{S , 1, v}^\perp$ for all $v \in S$ (see the discussion following Definition \ref{defring}), so $H^1_{\mc{D}_{S, 1}^\perp}(F, W_1^*) \to H^1_{\mc{D}_{S, ne}^{\perp}}(F, W_{ne}^*)[\vpi]$ is also an isomorphism. 
For any Taylor-Wiles prime $v$ of level $n$, $H^1(\gal{F_v}, W_1^*) \to H^1(\gal{F_v}, W_{ne}^*)[\vpi]$ is an isomorphism (since $(W_{ne}^*)^{\gal{F_v}} \to (W_{ne-1}^*)^{\gal{F_v}}$ is surjective in this case). Thus (c.f. the diagram on \cite[p. 522]{wiles:fermat}) if we can show there is a Taylor-Wiles set $Q_n$ of level $n$ such that $H^1(\gal{F, S}, W_1^*) \to \prod_{v \in Q_n} H^1(\gal{F_v}, W_1^*)$ is injective, then the same conclusion holds with $W_{ne}^*$ coefficients. For $W_1^*$, this is standard, and the argument of \cite[Theorem 2.49]{ddt} (see also \cite[Proposition 3.2.5]{kisin:bt}) shows we can always find such $Q_n$ with $|Q_n|=r:= h^1(\gal{F, S}, W_1^*)$ (this relies on Assumption \ref{twassumption}). In particular, for such $Q_n$, we have $h^1_{\mc{D}_{S \cup Q_n}^\perp}(\gal{F}, W_{ne}^*)=0$. 

With a bit more care, we can decrease the value of $r=|Q_n|$ as follows. If the preimage of $\mc{D}_{S , ne, v}^{\perp}$ under $W_1^* \to W_{ne}^*$ were $\mc{D}_{S , 1, v}^\perp$ for all $v \in S$, we could deduce that $H^1_{\mc{D}_{S,1}^\perp}(\gal{F, S}, W_1^*) \to H^1_{\mc{D}_{S, ne}^{\perp}}(F, W_{ne}^*)[\vpi]$ is an isomorphism, and then as above the argument of \cite[Theorem 2.49]{ddt} would apply to annihilate the dual Selmer group for $W_{ne}^*$ with $r= h^1_{\mc{D}_S^\perp}(\gal{F, S}, W_1^*)$ primes. This statement about the preimage is not, however, true in general. If we write $f \colon W_{ne} \to W_1$ and $f^* \colon W_1^* \to W_{ne}^*$ for the (dual) multiplication by $\vpi^{ne-1}$ maps (identifications using the choice of $\vpi$ left implicit), then a local class $\psi_v \in H^1(\gal{F_v}, W_1^*)$ satisfies $(f^* \psi_v) \in (\mc{D}_{ne, v})^\perp$ if and only if $\psi_v \in (f(\mc{D}_{ne,v}))^\perp$. Since each of our local conditions $\mc{D}_{m, v}$ is given as (the image in cohomology of) $\Hom(C_v, \mc{O}/\vpi^m)$ ($C_v$ as in the proof of Lemma \ref{gw}), $\dim(\mc{D}_{1, v}/f(\mc{D}_{ne, v}))$ is for $n \gg 0$ exactly the number of torsion factors $t_v$ appearing in Equation \eqref{localtorsion}. Thus a global class $\psi \in H^1_{\mc{D}_{S, ne}^{\perp}}(F, W_{ne}^*)[\vpi]$ belongs to a Selmer group $H^1_{L'_S}(\gal{F, S}, W_1^*)$ where for all $v \in S$, $\dim(L'_v/\mc{D}_{S,v, 1}^\perp)=t_v$. In particular, we can take 
\[
r= h^1_{\mc{D}_{S, 1}^\perp}(\gal{F, S}, W_1^*)+ \sum_{v \in S} t_v,
\]
as asserted in Theorem \ref{introthm}.

Whichever $r$ we take from above, combining Corollary~\ref{galcohgrowth} with Lemma~\ref{gw}, we find
\[
\begin{aligned}
h^1_{\mc D_S}(F,W_{ne})
&\leq h^1_{\mc D_{S\cup Q_n}}(F,W_{ne})-ner+d_r\\
&=h^1_{\mc D_{S\cup Q_n}^{\perp}}(F,W_{ne}^*)
  +c_{S,n}+d_r\\
&=c_{S,n}+d_r\\
&\leq c_S+d_r.
\end{aligned}
\]
%but for the local condition $\mc{D}_S$, we obtain 
%\[
%h^1_{\mc{D}_S}(\gal{F, S}, W_{ne})-h^1_{\mc{D}_S^\perp}(\gal{F, S}, W_{ne}^*)= c_{S, n},
%\]
%for the same constant $c_{S, n}$ appearing in Lemma \ref{gw}. Subtracting the two formulas, we find
%\[
%h^1_{\mc{D}_{S \cup Q_n}}(\gal{F, S \cup Q_n}, W_{ne})- h^1_{\mc{D}_S}(\gal{F, S}, W_{ne})= n e r- h^1_{\mc{D}_S^\perp}(\gal{F, S}, W_{ne}^*).
%\]

\end{proof}
\begin{cor}
Choose $r$ as in Proposition \ref{mainbound}. Then $H^1_{\mc{D}_S}(F, W_{E/\mc{O}})$ is finite with order bounded by $|k|^{c_S+d_r}$. Consequently, $H^1_g(\gal{F, S}, W_E)=0$.
\end{cor}

\begin{proof}
The first assertion follows immediately from Proposition
\ref{mainbound} upon taking the direct limit.  For the second, taking
the inverse limit over the cofinal sequence $m=ne$, the same bound gives
\[
\left(\varprojlim_m H^1_{\mc D_S}(F,W_m)\right)
 \otimes_{\mc O}E=0.
\]
By deformation theory, the left-hand side is naturally identified with
\[
\Hom_{\mc O}(\mf p/\mf p^2,E),
\]
where $\mf p$ is the kernel of the augmentation of
$R^{\mr{pst}}_{\br,S}$ corresponding to $\rho$.  This space parametrizes
self-extensions of $\rho\otimes_{\mc O}E$, with fixed determinant and
unramified outside $S$, whose restrictions at $v\mid p$ contain lattices
that are inverse limits of objects satisfying the torsion potentially
semistable condition used above.  As noted in the proof of Lemma \ref{gw}, Liu's proof of Fontaine's conjecture
\cite{liu:torsion}, applied after restriction to $F'_v$ and a Tate
twist, identifies this characteristic-zero local condition with
potential semistability over $F'_v$ with Hodge--Tate weights in $[a,b]$.
We also observed in Lemma~\ref{gw} that, for each $v\mid p$,
\cite[Lemma~1.2.5]{allen:selmer} identifies the image in
$H^1(F_v,W_E)$ of the tangent space of
$\overline R_v^\square[1/p]$ at the point defined by
$\rho|_{\gal{F_v}}$ with $H^1_g(F_v,W_E)$.  Since the local deformation
conditions at the places in $S$ away from $p$ are unrestricted, it
follows that
\[
\Hom_{\mc O}(\mf p/\mf p^2,E)
 \simeq H^1_g(\gal{F,S},W_E),
\]
and hence $H^1_g(\gal{F,S},W_E)=0$.
\end{proof}

We now return to the original field $F_0$ to deduce our final theorem.
\begin{thm}\label{finalthm}
Restriction induces an injection
\[
H^1_g(\gal{F_0},\ad^0(\rho_{\Pi_0,\iota})\otimes E)
\lhook\joinrel\longrightarrow
H^1_g(\gal{F},\ad^0(\rho_{\Pi,\iota})\otimes E).
\]
In particular, the Selmer group on the left vanishes. Moreover, if we define
local conditions $\mc{D}_{S_0}$ for deformations of
$\br_{\Pi_0, \iota}$ exactly as we did with $\mc{D}_S$ (see the discussion
after Definition \ref{defring}), then
\[
\lg H^1_{\mc{D}_{S_0}}(F_0, \ad^0(\rho_{\Pi_0, \iota}) \otimes E/\mc{O})
\leq c_S+d_r=[F:F_0]\sum_{v_0\in S_0}c_{v_0}+d_r,
\]
where $c_{v_0}$ is the local error term as defined in Lemma \ref{gw} (equal to $c_v$ for any $v \vert v_0$ since $F/F_0$ is split at $v_0$).
\end{thm}
\begin{proof}
The first assertion follows from inflation-restriction. With integral
coefficients, note that the local conditions $\mc{D}_{S_0,v_0}$ at
$v_0\in S_0$ restrict isomorphically to the local conditions
$\mc{D}_{S,v}$ for all $v\mid v_0$, since $F/F_0$ is split at $v_0$.
Since $[F:F_0]$ is either $1$ or $2$, and hence is prime to $p$, we see that for all $m\geq1$ the restriction 
\[
H^1_{\mc{D}_{S_0}}(F_0, \ad^0(\rho_{\Pi_0, \iota}) \otimes \mc{O}/\vpi^m) \xrightarrow{\sim} H^1_{\mc{D}_S}(F, W_m)^{\Gal(F/F_0)},
\] is an isomorphism. 
The length of the source is at most that of the target before taking
$\Gal(F/F_0)$-invariants. The result follows from Proposition
\ref{mainbound}.
%and the equality
%$c_S=[F:F_0]\sum_{v_0\in S_0}c_{v_0}$.
\end{proof}

\section{Implications and future directions}\label{generalsection}
We hope to generalize this work in two natural directions: to higher-rank $\Pi_0$, under usual Taylor-Wiles hypotheses (adequacy as in \cite{thorne:adequate}), and to more degenerate images $\br_{\Pi_0, \iota}$. We comment briefly on each of these.

If $\Pi_0$ is a regular algebraic polarizable (in the sense of \cite{blggt:potaut}) automorphic representation of $\mr{GL}_n(\mathbb{A}_F)$ for a CM field $F$, with $\br= \br_{\Pi_0, \iota}$ irreducible and $\br|_{\gal{F(\zeta_p)}}$ adequate, then we expect the template of the present paper to apply to prove vanishing of a suitable adjoint Bloch-Kato Selmer group (classifying polarized deformations). For this, one must combine the higher-rank arguments with Taylor-Wiles sets as in \cite{thorne:adequate} with estimates on congruence modules using the methods of \cite[\S 5, Appendix]{clozel-harris-taylor} (the desired lower bounds, as in part (5) of Proposition \ref{twsetup}, should not require Ihara's lemma). Note that the methods of the present paper require the weak multiplicity 1 condition from Lemma \ref{weakgenfree}, and in higher rank one makes a base-change (with some loss of integral information) to reduce the analysis to a unitary group where the relevant $\pi$ has multiplicity 1, using \cite[Th\'eor\`eme 5.4]{labesse:BC}. The output of this analysis, which we will report on soon, would be new quantitative bounds in higher rank for torsion adjoint Selmer groups, refining the vanishing after inverting $p$ first proven in this setting by Allen \cite{allen:selmer}.

Of course, even for $\mr{GL}_2$ over a totally real field, we have not
treated cases where $\br$ is reducible or where Taylor--Wiles sets as in
Proposition~\ref{mainbound} do not exist.  In this setting one must work with
sets $Q_n$ of ``approximate'' Taylor--Wiles primes: the Frobenius eigenvalues
need not be distinct modulo $\vpi$, but are distinct modulo $\vpi^M$ for some
$M$ depending on the image of $\rho_{\Pi,\iota}$ and independent of $n$.
One must also replace $R_{\br,S}^{\mr{pst}}$ by a suitable pseudodeformation
ring; Newton and Thorne developed the necessary techniques for cohomological control of these rings in \cite[\S~2]{newton-thorne:selmer}.

A significant obstacle is that, for approximate Taylor--Wiles sets, parts
(2)--(4) of Proposition~\ref{twsetup} fail.  Instead of an isomorphism
\[
M_{Q_n}/I_{Q_n}M_{Q_n}\simeq\mc{O}\oplus A_{\emptyset}
\]
with $A_{\emptyset}$ independent of $n$, one can probably only bound the
exponent of
$\bigl(M_{Q_n}/I_{Q_n}M_{Q_n}\bigr)_{\mr{tor}}$.  We have an approach to
handling this weaker control by modifying Propositions~\ref{twsetup} and
\ref{growth}, but it depends on the following hypothesis.

\begin{hypothesis}\label{CAconj}
Let $\overline I_{Q_n}$ be the image of $I_{Q_n}$ in
\[
\Hom_{\mc{O}}\bigl(M_{Q_n}/I_{Q_n}M_{Q_n},
I_{Q_n}M_{Q_n}/I_{Q_n}^2M_{Q_n}\bigr).
\]
Then
\[
\Fitt_{\mc{O}}(\overline I_{Q_n})
\subseteq\eta_{\lambda_{Q_n}}(M_{Q_n}).
\]
\end{hypothesis}

More abstractly, let $\lambda\colon T\to\mc{O}$ and $M$ be as in the
setup of Section~2, put $I=\ker(\lambda)$, and set
\[
\overline I_M:=
\im\left(
I/I^2\longrightarrow
\Hom_{\mc{O}}(M/IM,IM/I^2M)
\right).
\]
One may ask whether
\begin{equation}\label{general-CA-question}
\Fitt_{\mc{O}}(\overline I_M)\subseteq\eta_\lambda(M).
\end{equation}
For $M=T$, this is the standard elementary containment
\[
\Fitt_{\mc{O}}(I/I^2)\subseteq\eta_\lambda(T).
\]
In general, however, \eqref{general-CA-question} is false.\footnote{For
example, take $\mc O=\mathbb Z_3$,
\[
T=\mc O[x,y]/\bigl(x(x-3)(x-9)+3xy,\,
y(y-3)(y-9)+3xy\bigr),
\qquad \lambda(x)=\lambda(y)=0,
\]
and
\[
M=\mc O\left\langle
1,\frac{x+6}{9},
\frac{x^iy^j}{9}:0\leq i,j\leq2,\ 
(i,j)\notin\{(0,0),(1,0)\}
\right\rangle\subset T[1/3].
\]
Then $T$ is a reduced finite-flat local complete intersection of
$\mc O$-rank $9$, and $M$ is a $T$-stable lattice, generically free of
rank one. For $I=\ker(\lambda)$, one computes
$\overline I_M\simeq(\mc O/9)^2$ and
$\eta_\lambda(M)=(3^5)$. Hence
$\Fitt_{\mc O}(\overline I_M)=(3^4)\not\subseteq(3^5)
=\eta_\lambda(M)$. The module $M$ is not self-dual.  This example was found over several interactions with  ChatGPT.}

The automorphic modules $M_{Q_n}$ considered here are self-dual as
$\mathbb{T}_{Q_n}$-modules, and we know no counterexample to
\eqref{general-CA-question} for self-dual $M$.  Whether or not the precise
containment \eqref{general-CA-question} holds in that case, one may ask more
generally for lower bounds for $\lg(\overline I_M)$ in terms of invariants
intrinsic to $M$, such as its congruence ideal $\eta_\lambda(M)$.

If Hypothesis~\ref{CAconj} holds, or if a useful  weakening of it can
be established for the self-dual automorphic modules considered here, then
our methods might extend beyond the traditional Taylor--Wiles setting.
This could yield a new proof of the vanishing theorem of
\cite{newton-thorne:selmer}, together with explicit bounds on the
corresponding torsion Bloch--Kato Selmer groups.

\iffalse
\section{Almost Taylor-Wiles primes}
In this section we adapt the argument of \S \ref{twsection} to the much more
general setting where (in place of the assumption at the start of
\S \ref{twsection}) we only assume $\Pi$ is not induced from a quadratic
extension $K/F$ with $K$ contained in $F(\zeta_{p^\infty})$.
\fi

\section*{Acknowledgments}
We thank the referees for their helpful comments. C.B.K. thanks Sylvain Brochard and Srikanth Iyengar for many
discussions concerning the commutative algebra in the final section.

\section*{Statement on AI use}
C.B.K. used Chat GPT to generate the counterexample in \S \ref{generalsection}. We otherwise carried out the mathematical analysis and wrote the first version of this paper without use of AI. For final stage proofreading, we consulted a Chat GPT-generated ``referee report." 
\bibliographystyle{alpha}
\bibliography{biblio_selmer}

\end{document}